\documentclass[12pt]{amsart}

\newcommand{\Curve}{\mathcal{C}}

\newcommand{\Fdu}{\mathbf{F}}
\newcommand{\Fdup}{{\mathbf{F}_p}}
\newcommand{\Fduq}{{\mathbf{F}_q}}
\newcommand{\Fdupf}{{\mathbf{F}_{p^f}}}
\newcommand{\Fduqf}{{\mathbf{F}_{q^f}}}
\newcommand{\Fdur}{{\mathbf{F}_r}}
\newcommand{\Fduqn}{{\mathbf{F}_{q^n}}}
\newcommand{\Fdupbar}{{\overline{\mathbf{F}}_p}}

\newcommand{\Fdupbarn}{{\overline{\mathbf{F}}^n_p}}
\newcommand{\Fdupbarnone}{{\overline{\mathbf{F}}^{n+1}_p}}
\newcommand{\Zdu}{\mathbf{Z}}

\newcommand{\Qdu}{\mathbf{Q}}
\newcommand{\Qdubar}{\overline{\mathbf{Q}}}
\newcommand{\Rdu}{\mathbf{R}}
\newcommand{\Cdu}{\mathbf{C}}
\newcommand{\Qdul}{\mathbf{Q}_\ell}
\newcommand{\Qdulbar}{\overline{\mathbf{Q}}_\ell}
\renewcommand{\P}{\mathbf{P}}
\newcommand{\Adu}{\mathbf{A}}
\newcommand{\into}{\hookrightarrow}

\newcommand{\tensor}{\otimes}

\newcommand{\SLdu}{\mathrm{SL}}
\newcommand{\Spdu}{\mathrm{Sp}}
\newcommand{\GLdu}{\mathrm{GL}}
\renewcommand{\O}{\mathrm{O}}
\newcommand{\SOdu}{\mathrm{SO}}
\newcommand{\U}{\mathrm{U}}
\newcommand{\p}{{\mathfrak p}}
\newcommand{\qdu}{{\mathfrak q}}
\renewcommand{\d}{{\mathfrak d}}

\newcommand{\re}{{\rm Re\ }}

\DeclareMathOperator{\conddu}{Cond}
\DeclareMathOperator{\orddu}{ord}
\DeclareMathOperator{\rkdu}{Rank}

\DeclareMathOperator{\galdu}{Gal}
\DeclareMathOperator{\specdu}{Spec}

\DeclareMathOperator{\Frdu}{Fr}

\begin{document}

\title{Function fields and random matrices}
\author{Douglas Ulmer}
\maketitle

\footnotetext{April 20, 2006.  The author's research is partially supported by
  grants from the US National Science Foundation.}

\begin{quotation}
  ... le math\'ematicien qui \'etudie ces probl\`emes a l'impression
  de d\'echiffrer une inscription trilingue.  Dans la premi\`ere
  colonne se trouve la th\'eorie riemannienne des fonctions
  alg\'ebriques au sens classique.  La troisi\`eme colonne, c'est la
  th\'eorie arithm\'etique des nombres alg\'ebriques.  La colonne du
  milieu est celle dont la d\'ecouverte est la plus r\'ecente; elle
  contient la th\'eorie des fonctions alg\'ebriques sur un corps de
  Galois.  Ces textes sont l'unique source de nos connaissances sur
  les langues dans lesquels ils sont \'ecrits; de chaque colonne, nous
  n'avons bien entendu que des fragments; .... Nous savons qu'il y a
  des grandes diff\'erences de sens d'une colonne \`a l'autre, mais
  rien ne nous en avertit \`a l'avance.\hfil\break 
  \vspace{-0.15 in}\begin{flushright}A. Weil, ``De
    la m\'etaphysique aux math\'ematiques'' (1960)\end{flushright}
\end{quotation}

The goal of this survey is to give some insight into how
well-distributed sets of matrices in classical groups arise from
families of $L$-functions in the context of the middle column of
Weil's trilingual inscription, namely function fields of curves
over finite fields.  The exposition is informal and no proofs are
given; rather, our aim is to illustrate what is true by
considering key examples.

In the first section, we give the basic definitions and examples
of function fields over finite fields and the connection with
algebraic curves over function fields.  The language is a
throwback to Weil's Foundations, which is quite out of fashion but
which gives good insight with a minimum of baggage.  This part of
the article should be accessible to anyone with even a modest
acquaintance with the first and third columns of Weil's trilingual
inscription, namely algebraic functions on Riemann surfaces and
algebraic number fields.

The rest of the article requires somewhat more sophistication,
although not much specific technical knowledge.  In the second
section, we introduce $\zeta$- and $L$-functions over finite and
function fields and their spectral interpretation.  The cohomological
apparatus is treated purely as a ``black box.''  In the third section,
we discuss families of $L$-functions over function fields, the main
equidistribution theorems, and a small sample of applications to
arithmetic.  Although we do not give many details, we hope that this
overview will illuminate the function field side of the beautiful
Katz-Sarnak picture.

In the fourth section we give some pointers to the literature for
those readers who would like to learn more of the sophisticated
algebraic geometry needed to work in this area.

\section{Function fields}

In this first section we give a quick overview of function fields and
their connection with curves over finite fields.  The emphasis is on
notions especially pertinent to function fields over finite fields (as
opposed to function fields over algebraically closed fields), such as
rational prime divisors on curves, places of function fields, and
their behavior under extensions of fields and coverings of curves.
The section ends with a Cebotarev equidistribution theorem which is a
model for later more sophisticated equidistribution statements for
matrices in Lie groups.

\subsection{Finite fields}

If $p$ is a prime number, then $\Zdu/p\Zdu$ with the usual
operations of addition and multiplication modulo $p$ is a field
which we will also denote $\Fdup$.  If $\Fdu$ is a finite field,
then $\Fdu$ contains a subfield isomorphic to $\Fdup$ for a
uniquely determined $p$, the {\it
  characteristic\/} of $\Fdu$.  (The subfield $\Fdup=\Zdu/p\Zdu$ is the image
of the unique homomorphism of rings $\Zdu\to\Fdu$ sending 1 to 1.)
Since $\Fdu$ is a finite dimensional vector space over its
subfield $\Zdu/p\Zdu$, the cardinality of $\Fdu$ must be $p^f$ for
some positive integer $f$. Conversely, for each prime $p$ and
positive integer $f$, there is a field with $p^f$ elements, and
any two such are (non-canonically) isomorphic.  We may construct a
field with $q=p^f$ elements by taking the splitting field of the
polynomial $x^q-x$ over $\Fdup$.

It is old-fashioned but convenient to fix a giant field $\Omega$
of characteristic $p$ (say algebraically closed of infinite
transcendence degree over $\Fdup$) which will contain all fields
under discussion.  We won't mention $\Omega$ below, but all fields
of characteristic $p$ discussed are tacitly assumed to be
subfields of $\Omega$.  Given $\Omega$, we write $\Fdupbar$ for
the algebraic closure of $\Fdup$ in $\Omega$ (the set of elements of
$\Omega$ which are algebraic over $\Fdup$) and $\Fduq$ for the
unique subfield of $\Fdupbar$ with cardinality $q$.  Its elements
are precisely the $q$ distinct solutions of the equation
$x^q-x=0$. With this notation, $\Fduq\subset\Fdu_{q'}$ if and only
if $q'$ is a power of $q$, say $q'=q^k$ in which case $\Fdu_{q'}$
is a Galois extension of $\Fduq$ with Galois group cyclic of order
$k$ generated by the $q$-power Frobenius map $\Frdu_q(x)=x^q$.

\subsection{Function fields over finite fields}\label{ss:fields}
We fix a prime $p$.  A {\it function field $F$ of characteristic
$p$} is a finitely generated field extension of $\Fdup$ of
transcendence degree 1.  The {\it field of constants\/} of $F$ is
the algebraic closure of $\Fdup$ in $F$, i.e., the set of elements of
$F$ which are algebraic over $\Fdup$.  Since $F$ is finitely
generated, its field of constants is a finite field $\Fduq$.  When
we say ``$F$ is a function field over $\Fduq$'' we always mean
that $\Fduq$ is the field of constants of $F$. \vspace{0.2 in}

Examples:
\begin{enumerate}
\item The most basic example is the rational function field
$\Fduq(x)$
  where $q$ is a power of $p$ and $x$ is an indeterminate.  More
  explicitly, the elements of $\Fduq(x)$ are ratios of polynomials in
  $x$ with coefficients in $\Fduq$.  Its field of constants is $\Fduq$.
\item Let $q$ be a power of $p$, and let $F$ be the field extension of
  $\Fduq$ generated by two elements $x$ and $y$ and satisfying the
  relation $y^2=x^3-1$.  More precisely, let $F$ be the fraction field
  of $\Fduq[x,y]/(y^2-x^3+1)$ or equivalently $F=\Fduq(x)[y]/(y^2-x^3+1)$.
  The field of constants of $F$ is $\Fduq$.  If $p>3$, the field $F$ is
  not isomorphic to the rational function field $\Fduq(t)$.  (This is a
  fun exercise.  For hints, see [Sha77, p.~7].  Sadly, this point is
  missing from later editions of Shafarevitch's wonderful book.)  The
  cases $p=2$ and $p=3$ are degenerate: $F$ is isomorphic to the
  rational function field $\Fduq(t)$.  (If $p=2$, let $t=(y+1)/x$ and
  note that $x=t^2$ and $y=t^3-1$.  If $p=3$, let $t=y/(x-1)$ and note
  that $x=t^2+1$ and $y=t^3$.)
\item Similarly, if $p\neq2,5$ and $q$ is a power of $p$, let $F$ be
  the function field generated by $x$ and $y$ with relation
  $y^2=x^5-1$.  It can be shown that $F$ has field of constants $\Fduq$
  and is not isomorphic to either of the examples above.
\item Suppose that $p\equiv3\pmod4$ so that $-1$ is not a square in
  $\Fdup$.  Let $F$ be the function field generated over $\Fdup$ by
  elements $x_1,x_2,x_3$ with relations $x_1x_2=x_3$ and
  $x_2^2+x_3^2=0$.  It is not hard to see that the relations imply
  that $x_1^2=-1$ and so $F\cong \Fdu_{p^2}(x_2)$.  The moral is that
  the field of constants of $F$ is not always immediately visible from
  the defining generators and relations.
\end{enumerate}

If $F$ has constant field $\Fduq$, then any element $x\in
F\setminus\Fduq$ is transcendental over $\Fduq$ and so $F$
contains a subfield $\Fduq(x)$ isomorphic to the rational function
field. Since $F$ has transcendence degree 1, it is algebraic over
the subfield $\Fduq(x)$.

We can always choose the element $x\in F$ such that $F$ is a
finite {\it separable\/} extension of $\Fduq(x)$.  (It suffices to
choose $x$ which is not the $p$-th power of an element of $F$.)
The theorem of the primitive element then guarantees that $F$ is
generated over $\Fduq(x)$ by a single element $y$ satisfying a
separable polynomial over $\Fduq(x)$:
$$f(y)=y^n+a_1(x)y^{n-1}+\cdots+a_0(x)=0\qquad\text{with }a_i(x)\in\Fduq(x).$$
(Separable means that $f$ has distinct roots, or equivalently, $f$
and $\frac{df}{dy}$ are relatively prime in $\Fduq(x)[y]$.)  This
shows that $F$ is $\Fduq(x)[y]/(f(y))$.

More symmetrically, we may clear the denominators in the $a_i$ and
express the relation between $x$ and $y$ via a two-variable
polynomial over $\Fduq$:
$$g(x,y)=\sum b_{ij}x^iy^j=0\qquad\text{with }b_{ij}\in\Fduq.$$
This give us a presentation of $F$ as the fraction field of
$\Fduq[x,y]/(g(x,y))$.  Thus the general function field can be
generated over its constant field by two elements satisfying a
polynomial relation.  Note that this representation is far from
unique and it may be more natural in particular cases to give
several generators and relations.

\subsection{Curves over finite fields}\label{ss:curves}
Let $\Fdupbar$ be the algebraic closure of $\Fdup$ and let
$\P^n(\Fdupbar)$ denote the projective space of dimension $n$ over
$\Fdupbar$.  Thus elements of $\P^n(\Fdupbar)$ are by definition
the one-dimensional subspaces of the vector space $\Fdupbarnone$.  If
$(a_0,\dots,a_n)\in\Fdupbarnone\setminus(0,\dots,0)$, we write
$[a_0:\cdots:a_n]$ for the element of $\P^n(\Fdupbar)$ defined by
the subspace spanned by $(a_0,\dots,a_n)$.  We let $X_0,\dots,X_n$
denote the standard coordinates on $\Fdupbarnone$; of course the
$X_i$ do not give well-defined functions on $\P^n(\Fdupbar)$ but
the ratio of two homogenous polynomials in the $X_i$ of the same
degree gives a well-defined function on the set where the
denominator does not vanish.  In particular, on the subset
$X_0\neq0$, the functions $x_i=X_i/X_0$ ($i=1,\dots,n$) are a set
of coordinates which give a bijection between the set where
$X_0\neq0$ and the affine space $\Adu^n(\Fdupbar)=\Fdupbarn$.

We put a topology on $\P^n(\Fdupbar)$ by declaring that a
(Zariski) {\it
  closed\/} subset $Z\subset\P^n(\Fdupbar)$ is by definition the set of
points where some collection of homogeneous polynomials vanishes. We
may always take the set of polynomials to be finite and so a closed
set has the form
$$Z=\left\{[a_0:\cdots:a_n]\in\P^n(\Fdupbar)|
  f_1(a_0,\dots,a_n)=\cdots=f_k(a_0,\dots,a_n)=0\right\}$$
where $f_1,\dots,f_k\in\Fdupbar[X_0,\dots,X_n]$ are homogeneous
polynomials.  A closed subset $Z$ is said to be {\it defined over
  $\Fduq$\/} if we may take the $f_i$ to have coefficients in $\Fduq$.

We will work with the following definition, which is somewhat naive,
but suitable for our purposes: A (smooth, projective) {\it curve
  $\Curve$ over $\Fduq$\/} is a closed subset
$\Curve\subset\P^n(\Fdupbar)$ defined over $\Fduq$, such that:
\begin{enumerate}
\item $\Curve$ is infinite
\item there exist homogeneous polynomials $f_1,\dots,f_k$ vanishing
  identically on $\Curve$ such that for every $p\in\Curve$, the
  Jacobian matrix $(\frac{\partial f_i}{\partial X_j}(p))$
  ($i=1,\dots,k$ and $j=0,\dots,n$) has rank $n-1$
\item $\Curve$ is not the union of two proper closed subsets, i.e., if
  $Z_1$ and $Z_2$ are closed subsets and $\Curve=Z_1\cup Z_2$ then
  $\Curve= Z_1$ or $\Curve= Z_2$
\end{enumerate}
In the language of algebraic geometry, the first condition implies
that $\Curve$ has positive dimension and the first two conditions
imply that it is smooth and of dimension 1.  The third condition
says that $\Curve$ is {\it absolutely irreducible\/}.  If in the
third condition we insist that $Z_1$ and $Z_2$ be defined over
$\Fduq$ we arrive at the weaker condition that $\Curve$ is {\it
irreducible\/}. Although there are sometimes good reasons to
consider irreducible but not absolutely irreducible curves, for
simplicity we will not do so except in one example below.

We equip $\Curve$ with the Zariski topology induced from
$\P^n(\Fdupbar)$ so that its closed subsets are intersections of
$\Curve$ with closed subsets $Z\subset\P^n(\Fdupbar)$.

Warning: in the current literature a curve $\Curve$ is usually
defined in a more sophisticated way.  The set we are considering
here would be denoted $\Curve(\Fdupbar)$ and called the set of
$\Fdupbar$-valued points of $\Curve$. \vspace{0.2 in}

Examples:
\begin{enumerate}
\item $\P^1=\P^1(\Fdupbar)$ is the most basic example.  It is
defined by
  the zero polynomial on $\P^1$ (!) or, if that seems too
  tautological, by the equation $X_2=0$ in $\P^2(\Fdupbar)$.  Either
  representation makes it clear that $\P^1$ is defined over $\Fdup$.
\item[(1${}^{\prime}$)] For $p>2$, let $\Curve_2$ be the curve in
  $\P^2(\Fdupbar)$ defined over $\Fdup$ by the polynomial
  $X_1^2+X_2^2-X_0^2$.  Note that restricted to
  $\Curve_2\cap\{X_0\neq0\}$, the coordinate functions $x_i$ satisfy
  $x_1^2+x_2^2=1$.
\item[(1${}^{\prime\prime}$)] For any $p$, let $\Curve_3$ be the curve
  in $\P^3(\Fdupbar)$ defined over $\Fdup$ by the polynomials
  $X_0X_2-X_1^2$, $X_0X_3-X_1X_2$, and $X_1X_3=X_2^2$.  Note that
  restricted to $\Curve_3\cap\{X_0\neq0\}$, the coordinate functions
  $x_i$ satisfy $x_2=x_1^2$ and $x_3=x_1^3$.
\item Assume that $p>3$ and let $\Curve'_3$ be the curve in
  $\P^2(\Fdupbar)$ defined over $\Fdup$ by the polynomial
  $X_0X_2^2-X_1^3+X_0^3$.  (If $p=2$ or $3$ the second condition in
  the definition of a curve is not met: the Jacobian matrix is 0 at
  $[1:0:1]$ if $p=2$ and at $[1:1:0]$ if $p=3$.)  Note that restricted
  to $\Curve'_3\cap\{X_0\neq0\}$, the coordinate functions $x_i$
  satisfy $x_2^2=x_1^3-1$.
\item Assume that $p\neq2,5$ and let $\Curve_5$ be the closed subset
  of $\P^3(\Fdupbar)$ defined over $\Fdup$ by the equation polynomials
  $X_0X_2-X_1^2$, $X_0X_3^2-X_1X_2^2+X_0^3$, and
  $X_1X_3-X_2^3+X_0^2X_1$.  Note that restricted to
  $\Curve_5\cap\{X_0\neq0\}$, the coordinate functions $x_i$ satisfy
  $x_2=x_1^2$ and $x_3^2=x_1^5-1$.
\item Assume that $p\equiv3\pmod4$ so that $-1\in\Fdup$ is not a
square.
  Let $\Curve'_2$ be defined over $\Fdup$ by the three polynomials
  $X_0^2+X_1^2$, $X_2^2+X_3^2$, and $X_0X_3-X_1X_2$.  Then $\Curve'_2$
  is irreducible, but it is not absolutely irreducible and so it is
  not a curve by our definition.  Indeed, $\Curve'_2$ is the union of
  the two lines $\{X_0=iX_1, X_2=iX_3\}$ and $\{X_0=-iX_1,
  X_2=-iX_3\}$ defined over $\Fdu_{p^2}$ where $i\in\Fdu_{p^2}$ satisfies
  $i^2=-1$.  Note that restricted to $\Curve'_2\cap\{X_0\neq0\}$, the
  coordinate functions $x_i$ satisfy $x_1^2=-1$, $x_2^2+x_3^2=0$ and
  $x_3=x_1x_2$.
\end{enumerate}

\subsection{Morphisms and rational functions}\label{ss:maps}
If $\Curve\subset\P^m(\Fdupbar)$ and
$\Curve'\subset\P^{n}(\Fdupbar)$ are curves defined over $\Fduq$,
a {\it morphism\/} of curves is a map $\phi:\Curve\to\Curve'$ with
the property that at each point $P\in\Curve$, $\phi$ is
represented in an open neighborhood of $P$ by homogenous
polynomials.  In other words, for each $P\in\Curve$ there should
exist polynomials $f_0,\dots,f_{n}\in\Fdupbar[X_0,\dots,X_m]$, all
homogeneous of the same degree, such that for all $Q$ in some open
neighborhood of $P$, not all of the $f_i$ vanish at $Q$ and
$\phi(Q)=[f_0(Q):\cdots:f_{n}(Q)]$.  We say that $\phi$ is {\it
  defined over $\Fduq$\/} if it possible to choose the $f_i$ with
coefficients in $\Fduq$.  An {\it isomorphism\/} is a morphism
which is bijective and whose inverse is a morphism.
\vspace{0.2 in}

Examples:
\begin{enumerate}
\item If $f_0$ and $f_1$ are homogeneous polynomials in
$\Fduq[X_0,X_1]$
  of the same degree, not both 0, and with no common factors, then
$$[a_0:a_1]\mapsto[f_0(a_0,a_1):f_1(a_0,a_1)]$$
defines a morphism $\P^1\to\P^1$.  Using that $\Fduq[X_0,X_1]$ is
a unique factorization domain, one checks that every morphism from
$\P^1$ to itself defined over $\Fduq$ is of this form.
\item[(1${}^{\prime}$)] For $p>2$, the polynomials
$f_0=X_0^2+X_1^2$,
  $f_1=X_0^2-X_1^2$, and $f_2=2X_0X_1$ define a morphism from $\P^1$ over
  $\Fdup$ to the curve $\Curve_2$ in Example (1${}^{\prime}$) of
  Section~\ref{ss:curves}.  This morphism is an isomorphism with
  inverse defined by $f_0=\frac12(X_0+X_1)$ and
  $f_1=\frac12(X_0-X_1)$.
\item[(1${}^{\prime\prime}$)] For any $p$, the polynomials
  $f_0=X_0^3$, $f_1=X_0^2X_1$, $f_2=X_0X_1^2$, and $f_3=X_1^3$ define
  a morphism $\phi$ from $\P^1$ over $\Fdup$ to the curve $\Curve_3$ in
  Example (1${}^{\prime\prime}$) of Section~\ref{ss:curves}.  This
  morphism is an isomorphism with inverse defined on $\{X_0\neq0\}$ by
  $f_0=X_0$ and $f_1=X_1$ and on $\{X_3\neq0\}$ by $f_0=X_2$ and
  $f_1=X_3$.  In this example, it is not possible to define the
  inverse of $\phi$ by a single set of polynomials on all of
  $\Curve_3$.  Note also that the polynomials defining a morphism are
  in general not at all unique.  For example, on $\{X_0X_3\neq0\}$,
  the inverse of $\phi$ is defined both by $f_0=X_0$ and $f_1=X_1$ and
  by $f_0=X_2$ and $f_1=X_3$.
\item Let $\Curve'_3$ be as in Example (2) of
  Section~\ref{ss:curves}.  We define a morphism
  $\phi:\Curve'_3\to\P^1$ by setting $\phi([a_0:a_1:a_2])=[a_0:a_1]$
  on the open set where $a_0\neq0$ and
  $\phi([a_0:a_1:a_2])=[a_1^2:a_0^2+a_2^2]$ on the open set where
  $a_0^2+a_2^2\neq0$.  These requirements are compatible since if
  $a_0\neq0$ and $a_0^2+a_2^2\neq0$, then $a_1\neq0$ and
$$[a_0:a_1]=[a_0a_1^2:a_1^3]=[a_0a_1^2:a_0^3+a_0a_2^2]=[a_1^2:a_0^2+a_2^2].$$
If we think of $\P^1(\Fdupbar)\setminus\{[0:1]\}$ as $\Fdupbar$
via $[a_0:a_1]\mapsto a_1/a_0$, then the morphism $\phi$ extends
the function $x_1=X_1/X_0$, defined on $\Curve'_3\cap\{X_0\neq0\}$
to a morphism $\Curve'_3\to\P^1(\Fdupbar)$.  Again, it is not
possible to find a single pair of homogeneous polynomials
representing $\phi$ at all points of $\Curve'_3$.
\item[(2${}^{\prime}$)] Let $\Curve'_3$ be as above.  Choose a
  non-square element $a\in\Fdup$ ($p>3$) and define
  $\Curve^{\prime\prime}_3\subset\P^2(\Fdupbar)$ by the equation
  $aX_0X_2^2-X_1^3+X_0^3=0$.  Note that both $\Curve'_3$ and
  $\Curve^{\prime\prime}_3$ are defined over $\Fdup$.  Let
  $b\in\Fdu_{p^2}$ be a square root of $a$ and define a morphism
  $\phi:\Curve^{\prime\prime}_3\to\Curve'_3$ by
  $\phi([a_0:a_1:a_2])=[ba_o:a_1:a_2]$.  It is clear that $\phi$ is
  defined over $\Fdu_{p^2}$ and is an isomorphism.  On the other hand,
  one can show that $\Curve^{\prime\prime}_3$ and $\Curve_3'$ are not
  isomorphic over $\Fdup$.  This shows that two curves not isomorphic
  over their fields of definition may become isomorphic over a larger
  field.  One says that $\Curve^{\prime\prime}_3$ is a {\it twist\/}
  of $\Curve'_3$.
\item If $\Curve\subset\P^n(\Fdupbar)$ is a curve defined over
$\Fduq$,
  then there is an important morphism, the $q$-power Frobenius
  morphism $\Frdu_q:\Curve\to\Curve$, defined by
  $\Frdu_q([a_0:\cdots:a_n])=[a_0^q:\cdots:a_n^q]$.  Note that the
  fixed points of $\Frdu_q$ are precisely the points of $\Curve$ with
  coordinates in $\Fduq$.
\item If $\Curve\subset\P^n(\Fdupbar)$ is a curve and
$f_0,\dots,f_k$
  are homogeneous polynomials in $X_0,\dots X_n$ which do not all
  vanish identically on $\Curve$, then the map
  $\phi:\Curve\to\P^k(\Fdupbar)$ given by
 $$\phi([a_0:\cdots:a_n])=
[f_0(a_0,\dots,a_n):\cdots:f_k(a_0,\dots,a_n)]$$
is well-defined on the
  non-empty open subset of $\Curve$ where not all of the $f_i$ vanish.
  It is an important fact that $\phi$ can always be extended uniquely
  to a well-defined morphism on all of $\Curve$.  (NB: This is false
  for higher dimensional varieties.)  In particular, there are
  globally defined morphisms $x_i:\Curve\to\P^1$ extending the
  maps $[a_0:\cdots:a_n]\mapsto[a_0:a_i]$ which are {\it a priori\/}
  only defined on $\Curve\cap\{X_0\neq0\}$.
\end{enumerate}

A {\it rational function\/} on a curve $\Curve$ over $\Fduq$ is a
morphism $\phi:\Curve\to\P^1$ defined over $\Fduq$, except that we
rule out the constant morphism with image $\infty=[0:1]$.  (NB:
This is a reasonable definition only for curves, not for higher
dimensional varieties.)  In a neighborhood of any $P\in\Curve$,
$\phi$ can be represented by polynomials:
$\phi(Q)=[f_0(Q):f_1(Q)]$ where $f_0$ and $f_1$ are homogeneous of
the same degree and $f_0$ does not vanish identically.  It is
useful to think of $\phi$ as an $\Fdupbar$-valued function (with
poles) whose value at $Q$ is $\frac{f_1(Q)}{f_0(Q)}$. We say that
$\phi$ is {\it regular at $P\in\Curve$} if
$\phi(P)\neq\infty=[0:1]$.  If we restrict to an open set where
$\phi$ is regular, i.e., where $f_0$ does not vanish, then we get
a well-defined $\Fdupbar$-valued function.  If $\phi$ and $\phi'$
are two rational functions, we may restrict them to an open set
where they both give well-defined $\Fdupbar$-valued functions, add
or multiply them, and then extend back to rational functions on
$\Curve$.  More explicitly, if $\phi$ and $\phi'$ are represented
on some open set $U\subset\Curve$ by $[f_0:f_1]$ and
$[f_0':f_1']$, then $\phi+\phi'$ is represented by
$[f_0f'_0:f'_0f_1+f_0f'_1]$ and $\phi\phi'$ is represented by
$[f_0f'_0,f_1f'_1]$.  This gives the set of rational functions the
structure of a ring, in fact an algebra over $\Fduq$. This algebra
turns out to be a field extension of $\Fduq$ of transcendence
degree 1, i.e., a function field in the sense of the previous
subsection. It is denoted $\Fduq(\Curve)$.

Note that the ratio $f_1/f_0$ can be written as a rational
function (ratio of polynomials) in
$x_1=X_1/X_0,\dots,x_n=X_n/X_0$.  This shows that if
$\Curve\subset\P^n(\Fdupbar)$, then $\Fduq(\Curve)$ is generated
over $\Fduq$ by the rational functions $x_1,\dots,x_n$.  To
determine $\Fduq(\Curve)$, we need only determine the relations
among the $x_i$. \vspace{0.2 in}

Examples:
\begin{enumerate}
\item As noted above, a rational function on $\P^1$ is given by two
  homogeneous polynomials $f_0$ and $f_1$ of the same degree, with
  $f_0\neq0$.  Two rational functions $[f_0:f_1]$ and $[f_0':f_1']$
  are equal if and only if $f_1/f_0=f_1'/f_0'$.  Thus we see that
  rational functions on $\P^1$ are equivalent to rational functions
  (ratios of polynomials) in $x=X_1/X_0$, i.e., $\Fdup(\P^1)=\Fdup(x)$ and
  more generally $\Fduq(\P^1)=\Fduq(x)$.
\item[(1${}^{\prime}$)] The function fields of the curves $\Curve_2$
  and $\Curve_3$ in Examples (1${}^{\prime}$) and
  (1${}^{\prime\prime}$) of Section~\ref{ss:curves} are also
  isomorphic to $\Fdup(x)$.  One can see this by using the relations
  among the $x_i$ noted above, or by using the fact (to be explained
  below) that isomorphic curves have isomorphic function fields.
\item Let $\Curve'_3$ be as in Example (2) of Section~\ref{ss:curves}
  and let $x_1$ be the rational function $\phi$ of that example (so
  $x_1([a_0:a_1:a_2])=[a_0:a_1]$ or $[a_1^2:a_0^2+a_2^2]$).  Let $x_2$
  be the rational function defined on all of $\Curve'_3$ by
  $x_2([a_0:a_1:a_2])=[a_0:a_2]$.  Then $x_1$ and $x_2$ generate the
  field of rational functions on $\Curve'_3$ over $\Fduq$ and they
  satisfy the relation $x_2^2=x_1^3-1$.  In other words,
  $\Fduq(\Curve'_3)$ is the field in Example (2) of
  Section~\ref{ss:fields}.
\item Let $\Curve_5$ be as in Example (3) of Section~\ref{ss:curves}
  and define rational functions $x_1$ and $x_3$ by
$$x_1([a_0:a_1:a_2:a_3])=
\begin{cases}[a_0:a_1]&\text{if $a_0\neq0$}\cr
  [a_2^2:a_0^2+a_3^2]&\text{if $a_0^2+a_3^2\neq0$}\end{cases}$$ and
$$x_3([a_0:a_1:a_2:a_3])=[a_0:a_3].$$
(We leave it to the reader to check that these formulas do indeed
define rational functions on $\Curve_5$.)  It is not hard to see
that $x_1$ and $x_3$ generate $\Fduq(\Curve_5)$.  The equations
defining $\Curve_5$ imply that $x_3^2=x_1^5-1$ and that all
relations among $x_1$ and $x_3$ are consequences of this one. Thus
$\Fduq(\Curve_5)$ is the function field of Example (3) of
Section~\ref{ss:fields}.
\end{enumerate}

\subsection{The function field/curve dictionary}
The examples at the end of the last section illustrate the general
fact that if $\Curve$ is a curve defined over $\Fduq$, then the
field of rational functions $\Fduq(\Curve)$ is a function field,
i.e., a finitely generated extension of $\Fdup$ of transcendence
degree one, with field of constants $\Fduq$.

Conversely, it turns out that every function field $F$ with field
of constants $\Fduq$ is the field of rational functions of a curve
defined over $\Fduq$ which is uniquely determined up to
$\Fduq$-isomorphism.  We sketch one construction of the curve
corresponding to a function field $F$.  As we pointed out above,
$F$ may be generated over $\Fduq$ by two elements $x$ and $y$
satisfying a single relation
$$0=g(x,y)=\sum b_{ij}x^iy^j\qquad\text{with }b_{ij}\in\Fduq.$$
If $g$ has degree $d$, we form
$$G(X_0,X_1,X_2)=X_0^dg(X_1/X_0,X_2/X_0)=\sum b_{ij}X_0^{d-i-j}X_1^iX_2^j$$
and consider the closed subset of $\P^2(\Fdupbar)$ defined by
$G=0$. This closed subset will be infinite and irreducible, but it
will not in general be a curve under our definition, since it may
not satisfy the Jacobian condition.  If it does, we are finished.
If not, the closed set $\{G=0\}$ has singularities and the
classical process of blowing up (see [Ful89, Chap.~7]) gives an
algorithm to resolve the singularities and find a smooth curve in
some high-dimensional projective space with function field $F$. By
a suitable projection, the curve $\Curve$ can be embedded in
$\P^3(\Fdupbar)$.
In general we will not be able to find a {\it plane\/} curve with
function field $F$.  This is the case for example with the
function field in Example (3) of Section~\ref{ss:fields}
generated by $x$ and $y$ satisfying $y^2=x^5-1$.  The simplest
curve with this function field is a curve in $\P^3(\Fdupbar)$
defined by three equations.  

The dictionary between curves and function fields extends to
morphisms and field extensions.  More precisely, if $\Curve$ and
$\Curve'$ are two curves defined over $\Fduq$ and
$\phi:\Curve\to\Curve'$ is a non-constant morphism defined over
$\Fduq$, then composition with $\phi$ induces a ``pull-back''
homomorphism of fields $\Fduq(\Curve')\into\Fduq(\Curve)$ which is
the identity on $\Fduq$. Conversely, it can be shown that if $F$
and $F'$ are function fields over $\Fduq$ with corresponding
curves $\Curve$ and $\Curve'$, then a field inclusion $F'\into F$
which is the identity on $\Fduq$ is induced by a unique
non-constant morphism of curves $\Curve\to\Curve'$ which is
defined over $\Fduq$. \vspace{0.2 in}

Examples:
\begin{enumerate}
\item If $\Curve$ is a curve over $\Fduq$ and $x$ is a
non-constant
  rational function on $\Curve$, then $x$ is transcendental over
  $\Fduq$.  Thus the rational function field $F'=\Fduq(x)$ is a subfield
  of $F=\Fduq(\Curve)$.  The corresponding morphism $\Curve\to\P^1$ is
  the morphism $x$.
\item Suppose $F'$ is a function field with field of constants
$\Fduq$
  and $\Curve'$ is the corresponding curve over $\Fduq$.  If $r$ is a
  power of $q$ so that $\Fdur$ is a finite extension of $\Fduq$, then the
  function field $F=\Fdur F'$ corresponds to the same curve $\Curve'$
  viewed over $\Fdur$.  (Here $\Fdur F'$ is the compositum of $\Fdur$
  and $F'$, i.e., the smallest field containing both $\Fdur$ and
  $F'$.)  In other words, $F=\Fdur(\Curve')$.
\item We say that an extension of function fields $F/F'$ is {\it
    geometric} if it is separable and if the field of constants of $F$
  and $F'$ are the same.  If $n=[F:F']$ is the degree of the field
  extension, then the corresponding morphism of curves
  $\phi:\Curve\to\Curve'$ has degree $n$ in the sense that for all but
  finitely many $P\in\Curve'$, $\phi^{-1}(P)$ consists of $n$ points.
\item If $F/F'$ is a purely inseparable extension of function fields,
  say of degree $p^m$, then $F'=F^{p^m}$, the subfield of $p^m$-th
  powers.  In terms of suitable equations, the morphism
  $\Curve\to\Curve'$ acts on points by raising their coordinates to
  the $p^m$-th power.
\end{enumerate}
An arbitrary extension can be factored into three like these:
Given $F/F'$, let $\Fdur$ be the field of constants of $F$ and let
$F^{sep}$ be the separable closure of $F'$ in $F$.  Then $\Fdur
F'/F'$ is a constant field extension, $F^{sep}/\Fdur F'$ is
geometric, and $F/F^{sep}$ is purely inseparable.

\subsection{Points, prime divisors, and places}
As we have defined it, a curve $\Curve$ over $\Fduq$ is a set of
points with coordinates in $\Fdupbar$.  We would like to have a
set which reflects the fact that the equations defining $\Curve$
have coefficients in $\Fduq$.  The naive thing to look at would be
the set of $\Fduq$-rational points of $\Curve$, i.e., those with
coordinates in $\Fduq$, but this set is too small to be
useful---it may even be empty. The classical approach is to
consider $\Fduq$-rational prime divisors.

A {\it divisor\/} on $\Curve$ is a finite, formal, linear combination
$\d=\sum a_PP$ of points of $\Curve$ with integer coefficients.  A
divisor $\d$ is called {\it effective\/} if $a_P\ge0$ for all $P$.
The {\it degree\/} of $\d$ is $\deg(\d)=\sum a_P$.  The {\it
  support\/} of $\d$, written $|\d|$, is the set of points appearing
in $\d$ with non-zero coefficient.

If $\sigma\in\galdu(\Fdupbar/\Fduq)$ and $P\in\Curve$, then
$P^\sigma$ is again in $\Curve$.  (Here $\sigma$ acts on the
coordinates of $P$ and the claim follows from the fact that the
equations defining $\Curve$ have coefficients in $\Fduq$.)  We
extend this action to divisors by linearity ($(\sum
a_PP)^\sigma=\sum a_PP^\sigma$) and we say that a divisor $\d=\sum
a_PP$ is {\it $\Fduq$-rational\/} if it is fixed by the Galois
group, i.e., if $\d^\sigma=\d$ for all
$\sigma\in\galdu(\Fdupbar/\Fduq)$.

A {\it prime divisor\/} is an effective $\Fduq$-rational divisor
which is non-zero and cannot be written as the sum of two non-zero
$\Fduq$-rational effective divisors.  (Note that whether or not a
divisor is prime depends on the ground field over which we are
considering our curve.  A better terminology might be
$\Fduq$-prime, but we will stick with the traditional
terminology.)  It is not hard to see that the prime divisors of
$\Curve$ are in bijection with the orbits of
$\galdu(\Fdupbar/\Fduq)$ acting on $\Curve$.  If $\p$ is a prime
divisor, we define the residue field of $\p$ to be the field
generated over $\Fduq$ by the coordinates of any point in the
support of $\p$.  If $\p$ is prime and has degree $d$, then the
residue field at $\p$ is $\Fdu_{q^d}$.

If $P$ is a point of $\Curve$ and $f\in\Fduq(\Curve)$ is a
rational function on $\Curve$, then $f$ has a well-defined order
of vanishing or pole at $P$.  One motivation for considering prime
divisors is that the order of $f$ at $P$ is the same for all
points $P$ in the support of the prime divisor $\p$ containing
$P$.  In other words, the various points in $|\p|$ cannot be
distinguished from one another by the vanishing of
$\Fduq$-rational functions. \vspace{0.2 in}

Examples:
\begin{enumerate}
\item Let $\Curve=\P^1$ over $\Fduq$.  The divisors of degree 1
are
  simply the points of $\P^1$ with coordinates in $\Fduq$.  The prime
  divisors of degree $d>1$ are in bijection with the irreducible,
  monic polynomials in $\Fduq[x]$ of degree $d$, a polynomial
  corresponding to the formal sum of its roots.
\item Let $\Curve_3'$ be the curve in Example (2) of
  Section~\ref{ss:curves} over $\Fduq$.  If $a\in\Fduq$ with $a^3-1\neq0$,
  consider the points $P=[a:b:1]$ and $Q=[a:-b:1]$ where $b\in\Fdupbar$
  satisfies $b^2=a$.  The divisor $\d=P+Q$ has degree two and it is
  prime if and only if $b\not\in\Fduq$.  If $b\in\Fduq$, then $\d$ is the
  sum of two prime divisors, namely $P$ and $Q$.
\end{enumerate}

The set of prime divisors on $\Curve$ is more ``arithmetical''
than the full set of points on $\Curve$ (since it takes into
account that $\Curve$ is defined over $\Fduq$) and more convenient
and flexible than the set of $\Fduq$-rational points of $\Curve$.

Prime divisors play the role of the prime ideals of a number
field. More precisely, if $\p=\sum P_i$ is a prime divisor and if
$f\in\Fduq(\Curve)$ we say $f$ is regular (resp.~vanishes) at $\p$
if it is regular (resp.~vanishes) at one and therefore all of the
$P_i\in|\p|$.  The set of $f\in\Fduq(\Curve)$ which are regular at
$\p$ is a discrete valuation ring $R_\p$ with fraction field
$\Fduq(\Curve)$. The maximal ideal of $R_\p$ is the set of $f$
which vanish at $\p$. The residue field at $\p$ as we defined it
above turns out to be $R_\p$ modulo its maximal ideal.  We get a
valuation $\orddu_\p:\Fduq(\Curve)^\times\to\Zdu$ in the usual
way. It turns out that every non-trivial valuation of
$\Fduq(\Curve)$ is $\orddu_\p$ for a uniquely determined prime
divisor $\p$. (Therefore, it is possible, although not in my
opinion advisable, to eliminate the geometry completely and study
function fields via their valuations.  What one gains in algebraic
purity hardly seems to compensate for the loss of geometric
intuition this approach entails.)

Here is one respect in which the analogy between function fields
and number fields breaks down (``il y a des grandes diff\'erences
de sens d'une colonne \`a l'autre''): in a number field $F$, there
is a canonical Dedekind domain contained in $F$ whose primes give
the non-archimedean valuations of $F$, namely the ring of
integers.  In a function field, to get a Dedekind domain we fix a
non-empty set of prime divisors $S$ and then consider the ring $R$
of functions regular at all primes not in $S$.  The prime ideals
of $R$ are then in bijection with the prime divisors of
$\Fduq(\Curve)$ except those in $S$, and with the valuations of
$\Fduq(\Curve)$ except those arising from primes in $S$.  One
thinks of the primes in $S$ as the ``infinite primes'', but there
is no canonical choice for the set $S$.

\subsection{The Riemann-Roch theorem}
The Riemann-Roch theorem is true for curves over non-algebraically
closed fields and the statement is essentially the same as for the
case of curves over algebraically closed fields.  We give the basics
in our context.

Let $\Curve$ be a curve defined over $\Fduq$ with function field
$F=\Fduq(\Curve)$.  For each $P\in\Curve$ and $0\neq f\in F$,
there is a well-defined order of vanishing or pole of $f$ at $P$,
denoted $\orddu_P(f)$.  The {\it divisor of\/} $f$ is defined as
the formal sum $(f)=\sum_P\orddu_P(f)$ which is in fact a finite
sum. It is not hard to see that $(f)$ is $\Fduq$-rational and a
basic results says that it has degree 0: $\sum_P\orddu_P(f)=0$.

If $\d$ is an $\Fduq$-rational divisor, we define the Riemann-Roch
space $L(\d)$ by
$$L(\d)=\{f\in F^\times| (f)+\d\text{ is effective}\}\cup\{0\}.$$
Roughly speaking, $L(\d)$ consists of rational functions whose
poles are at worst given by $\d$.  It is clear that $L(\d)$ is an
$\Fduq$ vector space which turns out to be finite dimensional.
Note that $L(\d)$ is obviously zero if $\d$ has negative degree.

The Riemann-Roch theorem in its most basic form is a formula that
often allows one to compute the dimension $l(\d)$ of $L(\d)$.  The
theorem says that there is a non-negative integer $g$, the {\it genus
  of\/} $\Curve$ and a divisor $\omega$ of degree $2g-2$ such that for
all divisors $\d$
$$l(\d)-l(\omega-\d)=\deg(\d)-g+1.$$
The divisor $\omega$ is not unique (if $\omega$ works, then so does
$\omega+(f)$ for any non-zero $f$)).  Despite this ambiguity, $\omega$
is called a {\it canonical divisor\/}.  It turns out that $\omega$ can
be calculated as the divisor of a rational 1-form (i.e., a 1-form
possibly with poles) on $\Curve$.

It follows immediately that $l(\d)\ge\deg(\d)-g+1$ with equality if
$\deg\d>2g-2$.  This gives a large supply of functions with controlled
poles.

As an example, note that on $\P^1$ over $\Fduq$ the Riemann-Roch
space $L(d\infty)$ is just the space of polynomials of degree $d$,
which has dimension $d+1$.  It follows that the genus of $\P^1$ is
0.  One can check that the genus of the curve in Example~(2) of
Section~\ref{ss:curves} is 1 and the genus of the curve in
Example~(3) is 2.

As another application, which we leave as a simple exercise, the
theorem implies a partial converse to the statement that $\P^1$
has genus zero: if $\Curve$ has genus zero and an $\Fduq$-rational
divisor of degree 1, then $\Curve$ is isomorphic to $\P^1$.  It
turns out that over a finite field $\Fduq$ every curve has an
$\Fduq$-rational divisor of degree one, so this partial converse
is in fact a complete converse.

The reader curious about what a number field analog of the
Riemann-Roch theorem might be should consult Weil's ``Basic Number
Theory,'' [Wei95, Chap.~VI].

\subsection{Extensions, coverings, and splitting}\label{ss:splitting}
Let $\Curve$ and $\Curve'$ be curves defined over $\Fduq$ and let
$\phi:\Curve\to\Curve'$ be a morphism of curves defined over
$\Fduq$. We say that $\phi$ has degree $n$ if
$n=[\Fduq(\Curve):\Fduq(\Curve')]$. Given a point $P$ in $\Curve$
or $\Curve'$ we write $\Fduq(P)$ for the field generated over
$\Fduq$ by the coordinates of $P$.  We define an inverse image
mapping on divisors.  If $P\in\Curve'$ and if set-theoretically
the inverse image of $P$ in $\Curve$ is $\{Q_1,\dots,Q_k\}$, then
we assign a multiplicity $e_i$ to each $Q_i$ by choosing a
rational function $f$ vanishing simply at $P$ and setting $e_i=$
the order of vanishing of the pull-back of $f$ at $Q_i$.  We then
define $\phi^{-1}(P)$ as $\sum e_iQ_i$ and extend to divisors by
linearity.  It turns out that if $\Fduq(\Curve)$ is separable over
$\Fduq(\Curve')$ (i.e., if we have a geometric extension of
function fields), then for all but finitely many $P$, all the
$e_i$ are 1, and in general for all $P$, $\sum e_i=n$.

If $\p$ is a prime divisor of $\Curve$, then we may decompose
$\phi^{-1}(\p)$ into a sum of prime divisors
$\qdu_1,\dots,\qdu_g$. For each $\qdu_i$ we may define the {\it
residue degree\/} $f_i$ as $\deg\qdu_i/\deg\p$ or equivalently,
the degree of the field extension $\Fduq(Q)/\Fduq(P)$ where $P$ is
any point in the support of $\p$ and $Q$ is any point over $P$ in
the support of $\qdu_i$. The {\it ramification
  index\/} $e_i$ is the $e_i$ defined above for any point $P$ in the
support of $\p$ and any point $Q$ over $P$ in the support of
$\qdu_i$. It is a basic fact that for all $\p$,
$\sum_{i=1}^ge_if_i=n$ where $n=[F:F']$. \vspace{0.2 in}

Examples:
\begin{enumerate}
\item Suppose that $p>3$, $q$ is a power of $p$, $F$ is the fraction
  field of $\Fduq[x,y]/(y^2-x^3+1)$, and $\Fdu'=\Fduq(x)$, so that the
  corresponding morphism of curves $\phi:\Curve\to\Curve'=\P^1$ is as
  in Example (2) in Section~\ref{ss:maps}.  Suppose that $\p$ is a
  prime divisor of degree one corresponding to a finite $\Fduq$-rational
  point $P$ with coordinate $x=a$.  If $a^3-1=0$, then $\phi^{-1}(\p)$
  is a single prime $\qdu$ with $e=2$ and $f=1$; we say $\p$ is
  ramified.  If $a^3-1$ is a non-zero square of $\Fduq$, then
  $\phi^{-1}(\p)$ consists of two primes $\qdu_1$ and $\qdu_2$, both with
  $e=1$ and $f=1$; we say that $\p$ splits.  Finally, if $a^3-1$ is a
  non-square in $\Fduq$, then $\phi^{-1}(\p)$ consists of one prime $\qdu$
  with $e=1$ and $f=2$; we say that $\p$ is inert.
\item With notation as in the last example, if $\p$ is a general
  prime, say $\p=\sum P_i$, then the behavior of $\phi$ over each of
  the points $P_i$ is the same (one ramified point, two points with
  the same field of coordinates as $P_i$, or two points with
  coordinates in a quadratic extension of $\Fduq(P_i)$) and so
  $\phi^{-1}(\p)=2\qdu$ with $\deg\qdu=\deg\p$ ($\p$ ramifies),
  $\phi^{-1}(\p)=\qdu_1+\qdu_2$ with $\deg\qdu_i=\deg\p$ ($\p$ splits), or
  $\phi^{-1}(\p)=\qdu$ with $\deg\qdu=2\deg\p$ ($\p$ is inert).
\item If $\Curve$ is defined over $\Fduq$ and $r=q^n$, then we may
  consider the splitting of $\Fduq$-rational prime divisors into
  $\Fdur$-rational prime divisors.  This splitting is determined purely
  in terms of degrees: an $\Fduq$-rational prime $\p$ of degree $d$
  splits into $\gcd(n,d)$ $\Fdur$-rational primes, each with $e=1$ and
  $f=n/\gcd(d,n)$.
\item If $\Curve\to\Curve'$ is a morphism of curves defined over
$\Fduq$
  and is purely inseparable of degree $p^m$, then every prime $\p$ of
  $\Curve'$ pulls back to a single prime $\qdu$ of $\Curve$ with $e=p^m$
  and $f=1$.
\end{enumerate}
In the case of a morphism $\Curve\to\Curve'$ corresponding to a
geometric extension $F/F'$ which is Galois, it is easy to see that for
a fixed prime $\p$ of $\Curve'$, the ramification and residue degrees
$e_i$ and $f_i$ are all the same, in other words, $\p$ splits into $g$
primes, all with ramification index $e$ and residue degree $f$, and we
have $efg=n=[F:F']$.  Only finitely many $\p$ have $e>1$ and one can
make very precise statements about the distribution of primes having
allowable values of $f$ and $g$.  See Section~\ref{ss:cebo} below.

\subsection{Frobenius elements}\label{ss:frob}
Let $F'$ be a function field with constant field $\Fduq$ and let
$F$ be a finite Galois extension of $F'$ with Galois group $G$;
for simplicity we assume the extension $F/F'$ is geometric, i.e.,
the field of constants of $F$ is $\Fduq$.  Let
$\phi:\Curve\to\Curve'$ be the corresponding morphism of curves
over $\Fduq$.  Fix a finite extension $\Fdur$ of $\Fduq$ and a
point of $P\in\Curve'$ rational over $\Fdur$.  We may view $P$ as
an $\Fdur$-rational prime divisor.  Suppose that $\p_1,\dots,\p_g$
are the $\Fdur$-rational primes of $\Curve$ over $P$, so that as
divisors $\phi^{-1}(P)=e\p_1+\cdots+e\p_g$ where $e$ is the
ramification index.  The Galois group $G$ acts (transitively in
fact) on the set of $\p_i$ and we let $D_{\p_i}\subset G$ denote
the stabilizer of $\p_i$, the {\it decomposition group at\/}
$\p_i$.  Then $D_{\p_i}$ acts on the residue field at $\p_i$ and
so we have a homomorphism $D_{\p_i}\to\galdu(\Fdu_{r'}/\Fdur)$
where $\Fdu_{r'}=\Fdur(\p_i)=\Fdur(Q)$ for any $Q\in|\p_i|$.  This
homomorphism is surjective with kernel denoted $I_{\p_i}$, the
{\it inertia group
  at\/} $\p_i$.  It turns out that the order of the inertia group is
$e$, the ramification index of $\p_i$.  When $e=1$, there is a
distinguished element of $D_{\p_i}$, namely the one that maps to
the $r$-power Frobenius in $\galdu(\Fdu_{r'}/\Fdur)$.  When $e>1$
we get a distinguished coset of $I_{\p_i}$ in $D_{\p_i}$.
Changing the choice of $\p_i$ changes $D_{\p_i}$, $I_{\p_i}$ and
the distinguished element or coset by conjugation by an element of
$G$.  Therefore, we get a well-defined conjugacy class in $G$
depending only on $\Fdur$ and $P$ which we denote $\Frdu_{\Fdur,P}$.
Similarly, we write $D_{\Fdur,P}$ and $I_{\Fdur,P}$ for the
conjugacy classes of subgroups of $G$ defined as above.  It is not
hard to check that $\Frdu_{\Fdu_{r^n},P}=\Frdu_{\Fdur,P}^n$.

One also associates decomposition and inertia subgroups and a
Frobenius element to a prime $\p$ of $\Curve$ as follows: we let
$\Fdur$ be the residue field at $\p$ and choose $P\in|\p|$ and
then set $D_\p=D_{\Fdur,P}$, $I_\p=I_{\Fdur,P}$, and
$\Frdu_\p=\Frdu_{\Fdur,P}$.  The resulting conjugacy classes are
well-defined independently of the choice of $P$.  This Frobenius
is more analogous to the Frobenius element considered over
number fields.

Example: Let $\Curve\to\Curve'=\P^1$ be the morphism considered in
Example (2) in Section~\ref{ss:maps} and again in Example (2) in
Section~\ref{ss:splitting}.  This is a Galois covering with group
$G=\Zdu/2\Zdu$.  If $a\in\Fdur$ is such that $a^3-1\neq0$, and
$P\in\P^1$ is the point $[1:a]$, then the Frobenius class $\Frdu_P$
is 1 if $a^3-1$ is a square in $\Fdur$ and is $-1$ if it is not a
square.  If $\p$ is an $\Fduq$-rational prime divisor of $\P^1$,
then $\Frdu_\p$ is 1 if $\p$ splits and is $-1$ if $\p$ is inert.

The definitions of decomposition and inertia subgroups and Frobenius
elements extend to infinite Galois extensions in exactly the same way
as in the number field context.

\subsection{Cebotarev equidistribution}\label{ss:cebo}
The classical Cebotarev density theorem says roughly that Frobenius
elements are equidistributed in the Galois group of a Galois extension
of number fields.  To discuss a function field analogue, we keep the
notations of the last section so that $F/F'$ is a geometric Galois
extension of function fields over $\Fduq$, with corresponding morphism
of curves $\Curve\to\Curve'$ defined over $\Fduq$.  We consider the
distribution of Frobenius conjugacy classes $\Frdu_{\Fdur,P}$ as $P$
varies over $\Fdur$-rational points of $\Curve'$ for large $r$.

One analogue of the Cebotarev density theorem for function fields says
that the Frobenius classes become equidistributed as $r$ tends to
infinity.  More precisely, if $C\subset G$ is a conjugacy class, then
\begin{equation*}
\lim_{r\to\infty}\frac{|\{P\in\Curve'(\Fdur)|\Frdu_{\Fdur,P}\in C\}|}
{|\{P\in\Curve'(\Fdur)\}|}=\frac{|C|}{|G|}
\end{equation*}
where $r$ tends to infinity through powers of $q$.  A useful way to
rephrase this is to consider conjugation invariant functions $f$ on
$G$.  It make sense to evaluate such a function on a Frobenius
conjugacy class and we have
\begin{equation*}
\lim_{r\to\infty}\left|\frac1{|\Curve'(\Fdur)|}\sum_{P\in\Curve'(\Fdur)}f(\Frdu_{\Fdur,P})
-\frac1{|G|}\sum_{g\in G}f(g)\right|=0
\end{equation*}

There is a more precise statement about the rate of convergence: given
data as above, there exists a constant depending only on $F/F'$ and
$f$ such that for all powers $r$ of $q$,
\begin{equation*}
\left|\frac1{|\Curve'(\Fdur)|}\sum_{P\in\Curve'(\Fdur)}f(\Frdu_{\Fdur,P})
-\frac1{|G|}\sum_{g\in G}f(g)\right|\le Cr^{-1/2}.
\end{equation*}
The constant $C$ can be made quite explicit in terms of the
representation theory of $G$ and the expansion of $f$ in terms of
characters.  See [KS99b, 9.7.11-13] for details.

As a very simple example of what this means in down-to-earth
terms, we return to Example (2) of Section~\ref{ss:splitting}.
In that context, Cebotarev equidistribution says that for large
$r$, for about $1/2$ of the elements $a\in\Fdur$, $a^3-1$ is a
square and for about $1/2$ of the $a$, it is not a square.

\section{$\zeta$-functions and $L$-functions}\label{s:zetas}

In this section we define $\zeta$- and $L$-functions, give some
examples, and discuss the spectral interpretation.  Warning: we use a
non-standard, radically simplified notation for certain cohomology
groups.  See Section~\ref{s:refs} for references with a more complete
treatment.

\subsection{The $\zeta$-function of a curve}
Let $F$ be a function field with field of constants $\Fduq$.  Let
$\Curve$ be the corresponding curve and denote by $\Curve^0$ the
set of $\Fduq$-rational prime divisors of $\Curve$.  We define the
zeta-function of $\Curve$ in analogy with the Riemann
zeta-function:
$$\zeta(\Curve,s)=\prod_{\p\in\Curve^0}(1-N\p^{-s})^{-1}$$
where $N\p=q^{\deg\p}$ is the number of elements in the residue field
at $\p$.
(This function depends not just on the curve $\Curve$ but also on
the constant field $\Fduq$ and when we want to make this
dependence explicit, we write $\zeta(\Curve/\Fduq,s)$.)

If $C_m$ denotes the number of primes in $\Curve^0$ of degree $m$
and $N_n$ denotes the number of points of $\Curve$ defined over
$\Fduqn$, then we have
$$N_n=\sum_{m|n}mC_m.$$
Rearranging formally, we find that
$$\zeta(\Curve,s)=\exp\left(\sum_{n=1}^\infty\frac{N_n}nq^{-ns}\right)$$
which makes the diophantine interest of $\zeta$ quite visible.

The product defining $\zeta(\Curve,s)$ and the rearranged sum
converge absolutely in the region $\re s>1$.  Using the
Riemann-Roch theorem, one can show that $\zeta(\Curve,s)$ extends
to a meromorphic function on all of $\Cdu$, with simple poles at
$s=1$ and $s=0$ and holomorphic elsewhere, and that it satisfies a
functional equation relating $s$ and $1-s$.  (There are no
$\Gamma$-factors because the product defining $\zeta$ is over all
places of $F$.)  More precisely,
$$q^{-s(1-g)}\zeta(\Curve,s)=q^{(s-1)(1-g)}\zeta(\Curve,1-s)$$
where $g$ is the genus of $\Curve$.

Here are some examples: If $F$ is the rational function field with
constant field $\Fduq$, so that $\Curve=\P^1$, then $N_n=q^n+1$
and so
$$\zeta(\Curve,s)=\frac{1}{(1-q^{-s})(1-q^{1-s})}.$$

Let $\Curve$ be the curve with affine equation $y^2=x^3-x$ over
$\Fdup$ where $p\equiv3\pmod4$ and $p>3$.  Using the fact that
$-1$ is not a square modulo $p$, it is easy to check that the
number of points on $E$ over $\Fdup$ is $p+1$ and more generally,
if $f$ is odd, the number of points on $E$ with coordinates in
$\Fdu_{p^f}$ is $p^f+1$.  (One considers pairs $x=a$ and $x=-a$,
excluding $x=0$ and $\infty$.  Since $-1$ is not a square in
$\Fduq$, $x^3-x$ is a square for exactly one of $x=a$ or $x=-a$;
when it is a square there are two values of $y$ with $y^2=x^3-x$
and none when it is not.  Thus the number of solutions with finite
non-zero $a$ is $q-1$ and the total number of solutions is $q+1$.)
A somewhat more elaborate argument using exponential sums allows
one to show that for even $f$, the number of solutions over
$\Fdupf$ is $p^f+1-2(-p)^{f/2}$.  (See Koblitz [Kob93, II.2] or
Ireland and Rosen [IR90, Chap.~18] for a nice exposition of this
argument.) Using the expression for $\zeta$ in terms of the $N_n$,
we conclude that
\begin{equation*}
\zeta(\Curve/\Fdup,s)=\frac{(1-\sqrt{-p}p^{-s})(1+\sqrt{-p}p^{-s})}
{(1-p^{-s})(1-p^{1-s})} =\frac{1+p^{1-2s}}{(1-p^{-s})(1-p^{1-s})}.
\end{equation*}

As a third example we assume that $p>2$ and $q=p^f\equiv1\pmod 3$ and
consider the curve $\Curve$ with affine equation $y^3=x^4-x^2$, or
rather the smooth, projective curve obtained from this one by
desingularization.  (This curve is singular at $(x,y)=(0,0)$, but
there is exactly one point over this one in the smooth curve, so for
the purposes of counting points we may ignore this.)  This curve has
genus $g=2$.

Let $\lambda:\Fduq^\times\to\Cdu^\times$ be a character of order
exactly 6 and for $a=1,2,4,5$ define
$$J_a=\sum_{\substack{x\in\Fduq\\x\neq0,1}}\lambda^a(x(1-x)).$$
It is not hard to check that $|J_i|=q^{1/2}$ and
$J_5=\overline{J}_1$, $J_4=\overline{J}_2$.  Using arguments
similar to those in Koblitz or Ireland and Rosen, one verifies
that the number of points on $\Curve$ over $\Fduqf$ is
$q^f+1-\sum_{a\in\{1,2,4,5\}}J_a^f$.  This implies that
$$\zeta(\Curve/\Fduq,s)
=\frac{\prod_{a\in\{1,2,4,5\}}(1-J_aq^{-s})}{(1-q^{-s})(1-q^{1-s})}.$$

In general, if $\Curve$ has genus $g$ then $\zeta(\Curve,s)$ has the
form
$$\frac{P(q^{-s})}{(1-q^{-s})(1-q^{1-s})}$$
where $P$ is a polynomial of degree $2g$ with integer coefficients
and constant term 1.  Writing
$P(T)=\prod_{i=1}^{2g}(1-\alpha_iT)$, the functional equation for
$\zeta$ is equivalent to the fact that the set of inverse roots
$\alpha_i$ is invariant under $\alpha_i\mapsto q/\alpha_i$.
Moreover, $\zeta$ satisfies an analogue of the Riemann hypothesis:
all of the inverse roots $\alpha_i$ have absolute value $q^{-1/2}$
and so the zeros of $\zeta$ lie on the line $\Re(s)=1/2$. These
results were proven in general by Weil in [Wei48].

More generally, one can define a zeta function for any variety defined
over a finite field via a product or exponentiated sum as above.  If
$X$ is smooth and complete of dimension $d$, then one knows that
$\zeta(X,s)$ is a rational function in $q^{-s}$ of a very special
form.  More precisely,
$$\zeta(X,s)=
\frac{P_1(q^{-s})P_3(q^{-s})\cdots
  P_{2d-1}(q^{-s})}{P_0(q^{-s})P_2(q^{-s})\cdots P_{2d}(q^{-s})}$$
where each $P_i$ is a polynomial with integer coefficients all of
whose inverse roots have complex absolute value $q^{i/2}$ (an analogue
of the Riemann hypothesis).  Moreover, if the inverse roots of $P_i$
are $\alpha_1,\dots,\alpha_k$, then the inverse roots of $P_{2d-i}$
are $q^d/\alpha_1,\dots,q^d/\alpha_k$ and so $\zeta(X,s)$ extends to a
meromorphic function in the plane and satisfies a functional equation
for $s\to d-s$.  These properties of the $\zeta$-function were
conjectured by Weil in [Wei49] and proved in full generality by Deligne
in 1974.

\subsection{Spectral interpretation of $\zeta$-functions}
Already at the time he made his famous conjectures, Weil envisioned a
cohomological explanation for the conjectured properties of the zeta
function.  This was provided in important cases by Weil and later in
full generality by Grothendieck, Deligne, and collaborators.

We fix an auxiliary prime $\ell$ not equal to the characteristic
of $\Fduq$.  Attached to a curve $\Curve$ over a finite field
$\Fduq$ are finite-dimensional $\Qdul$-vector spaces
$H^0(\Curve)$, $H^1(\Curve)$ and $H^2(\Curve)$ each equipped with
an action of $\galdu(\Fdupbar/\Fduq)$. The $\zeta$-function of
$\Curve$ then has an interpretation in terms of the spectrum of
the $q$-power Frobenius $\Frdu_q$, which is a generator of
$\galdu(\Fdupbar/\Fduq)$, namely
$$\zeta(\Curve,s)=\frac{P_1(q^{-s})}{P_0(q^{-s})P_2(q^{-s})}$$
where
$$P_i(T)=\det\left(1-T\Frdu_q|H^i(\Curve)\right).$$

(It turns out that the eigenvalues of $\Frdu_q$ are algebraic
numbers, so that we may interpret them as complex numbers.  In
fact the coefficients of the reversed characteristic polynomials
appearing here are integers, so there is no dependence on an
embeddings of $\Qdubar$ into $\Qdulbar$ and $\Cdu$.)

It turns out that $H^0(\Curve)$ is one-dimensional with trivial action
of $\Frdu_q$, $H^2(\Curve)$ is one-dimensional with $\Frdu_q$ acting by
multiplication by $q$ and $H^1(\Curve)$ is $2g$-dimensional, where $g$
is the genus of $g$.  This shows that $\zeta(\Curve,s)$ is a rational
function in $q^{-s}$ of the form mentioned in the last section.

The functional equation is a manifestation of a Poincar\'e duality:
there are pairings $H^i(\Curve)\times H^{2-i}(\Curve)\to H^2(\Curve)$
compatible with the actions of $\Frdu_q$ and this shows that the
eigenvalues of $\Frdu_q$ on $H^i$ are $q$ divided by the eigenvalues of
$\Frdu_q$ on $H^{2-i}$, which is the content of the functional equation.

The Riemann hypothesis, namely that the zeros of $\zeta(\Curve,s)$
lie on the line $\Re(s)=1/2$, is equivalent to the statement that
the eigenvalues of $\Frdu_q$ on $H^1(\Curve)$ have complex absolute
value $q^{1/2}$.

All of the above generalizes to smooth proper varieties of any
dimension over $\Fduq$.  For an $X$ of dimension $d$, there are
finite-dimensional $\Qdul$-vector spaces $H^0(X),\dots,H^{2d}(X)$
with an action of $\Frdu_q$; $H^0(X)$ is one-dimensional with trivial
$\Frdu_q$ action and $H^{2d}(X)$ is one-dimensional with $\Frdu_q$
acting by multiplication by $q^d$.  There is a Poincar\'e duality
pairing $H^i(X)\times H^{2d-i}(X)\to H^{2d}(X)$ which is
non-degenerate and compatible with the Frobenius actions. Finally,
the eigenvalues of $\Frdu_q$ on $H^i(X)$ are algebraic integers with
absolute value $q^{i/2}$ in every complex embedding.

\subsection{Examples of $L$-functions}\label{ss:Lexamples}
Just as in the number field case, we can define $L$-functions
associated to representations of the absolute Galois group of a
function field.  Before giving the general definitions, we consider
three examples.

First, let $F$ be a quadratic extension of $\Fduq(t)$,
corresponding to a branched cover $\Curve\to\P^1$ of degree 2.
Since $F/\Fduq(t)$ is a Galois extension with group $\{\pm1\}$, we
get a quadratic character
$$\chi:\galdu(\overline{\Fduq(t)}/\Fduq(t))\to\galdu(F/\Fduq(t))\to\{\pm1\}.$$
Let us define the $L$-function of $\chi$ as
$$L(\chi,s)=\prod_{\p\in(\P^1)^0}(1-\chi(\p)N\p^{-s})^{-1}$$
where for unramified $\p$, $\chi(\p)=\chi(\Frdu_\p)$ is 1 if $\p$ splits
in $F$ and $-1$ if $\p$ is inert; we set $\chi(\p)=0$ if $\p$ is
ramified in $F$.  An elementary (Euler-factor by Euler-factor)
computation shows that
$$\zeta(\Curve,s)=\zeta(\P^1,s)L(\chi,s).$$
On the other hand,
$$\zeta(\Curve,s)=\frac{P(q^{-s})}{(1-q^{-s})(1-q^{1-s})}$$
and
$$\zeta(\P^1,s)=\frac{1}{(1-q^{-s})(1-q^{1-s})}$$
and so
$$L(\chi,s)=P(q^{-s}).$$
The functional equation for $\zeta$ is equivalent to
$$q^{gs}L(\chi,s)=q^{g(1-s)}L(\chi,1-s).$$

This applies in particular to the curve $y^2=x^3-x$ considered above:
we view it as a degree two cover of the $t$-line by $(x,y)\mapsto
t=x$.  It follows that
$$L(\chi,s)=(1-\sqrt{-p}p^{-s})(1+\sqrt{-p}p^{-s})
=1+p^{1-2s}.$$

For a second class of examples, consider a Galois extension
$F/\Fduq(t)$ with Galois group $\Zdu/d\Zdu$, corresponding to a
degree $d$ cyclic covering of curves $\Curve\to\P^1$.  Let
$\chi:\galdu(\overline{\Fduq(t)}/\Fduq(t))\to
\galdu(F/\Fduq(t))\to\mu_d\subset\Qdubar^\times$ be a complex
valued character of order exactly $d$ and for $i=1,\dots,d-1$
define
$$L(\chi^i,s)=\prod_{\p\in(\P^1)^0}(1-\chi^i(\p)N\p^{-s})^{-1}$$
where $\chi(\p)=\chi(\Frdu_\p)$ for unramified $\p$ and $\chi(\p)=0$ if
$\p$ is ramified in $F$.  Again an elementary calculation shows that
$$\zeta(\Curve,s)=\zeta(\P^1,s)L(\chi,s)L(\chi^2,s)\cdots L(\chi^{d-1},s).$$
It turns out that each $L(\chi^i,s)$ for $i=1,\dots,d-1$ is a
polynomial in $q^{-s}$ and their product is the numerator $P(q^{-s})$
of $\zeta(\Curve,s)$.

For $d>2$ a new phenomenon becomes apparent: the functional equation
links two distinct $L$-functions.  More precisely, we have
$$q^{N_is/2}L(\chi^i,s)=\epsilon q^{N_i(1-s)/2}L(\chi^{-i},1-s)$$
where $N_i=N_{-i}$ is the degree of $L(\chi^i,s)$ as a polynomial in
$q^{-s}$ and $\epsilon$ is a complex number of absolute value 1.  This
will be important later when we discuss symmetry types.

As a specific example of this type, we consider the curve $\Curve$
defined by $y^3=x^4-x^2$, discussed above, viewed as a Galois
cover of $\P^1$ of degree 3 via $(x,y)\mapsto t=x$.  For a
suitable choice of character $\chi:\galdu(F/\Fduq(t))\to\mu_3$, we
have $L(\chi,s)=(1-J_1q^{-s})(1-J_4q^{-s})$ and
$L(\chi^2,s)=(1-J_2q^{-s})(1-J_5q^{-s})$.

A third, more elaborate, class of examples comes from elliptic curves.
Let $E$ be an elliptic curve defined over $F$.  This could be given,
for example, by a Weierstrass equation
$$y^2+a_1xy+a_3y=x^3+a_2x^2+a_4x+a_6$$
where the $a_i$ are in $F$.  If $E$ has good reduction at a place $\p$
of $F$, we define a local Euler factor by
$$L_\p(E,s)=(1-a_\p q_\p^{-s}+q_\p^{1-2s})$$
where $q_\p$ is the cardinality of the residue field at $\p$ and
$q_\p-a_\p+1$ is the number of points on the reduction of $E$ at $\p$.
If $E$ has bad reduction at $\p$, we define a local factor by
$$L_\p(E,s)=\begin{cases}
1-q_\p^{-s}&\text{if $E$ has split multiplicative reduction at $\p$}\\
1+q_\p^{-s}&\text{if $E$ has non-split multiplicative reduction at $\p$}\\
1&\text{if $E$ has additive reduction at $\p$.}
\end{cases}
$$
Then we define the global (Hasse-Weil) $L$-function of $E$ as
$$L(E,s)=\prod_{\p}L_\p(E,s)^{-1}.$$
This $L$-function turns out to be a rational function in $q^{-s}$
and it satisfies a functional equation for $s\to2-s$.  More
precisely, if $E$ is not isomorphic to an elliptic curve defined
over $\Fduq$, then $L(E,s)$ is a polynomial in $q^{-s}$ whose
degree is determined by the genus of the curve corresponding to
$F$ and the places of bad reduction of $E$.  In this case,
$$L(E,s)=\prod_{i=1}^N(1-\alpha_iq^{-s})$$
where the set of inverse roots $\alpha_i$ is invariant under
$\alpha_i\mapsto q^2/\alpha_i$ and each of them has complex
absolute value $q$.  In particular, the zeros of $L(E,s)$ lie on
the line $\Re(s)=1$.

\subsection{$L$-functions attached to Galois representations}
As in the number field case, over function fields there are two
general classes of $L$-functions, automorphic $L$-functions attached
to automorphic representations (generalizing Dirichlet characters,
Hecke characters, etc.) and ``motivic'' $L$-functions attached to
representations of Galois groups, and a Langlands philosophy which
very roughly speaking says that the latter are the same as the former.
In the function field setting there is a quite satisfactory
understanding of the analytic properties of motivic $L$-functions
which we sketch in this and the following section.

As usual, let $F=\Fduq(\Curve)$ be the function field of a curve
over $\Fduq$.  We fix a prime $\ell$ and write $E$ for a finite
extension of $\Qdul$ which we may expand as necessary in the
course of the discussion.  The basic input data is a
representation
$$\rho:\galdu(\overline{F}/F)\to\GLdu_n(E)$$
which is continuous (for the Krull topology on
$\galdu(\overline{F}/F)$ and the $\ell$-adic topology on
$\GLdu_n(E)$) and unramified outside a finite set of places of
$F$.  The latter means that for all but finitely many primes $\p$,
$\rho(I_\p)=\{1\}$ where $I_\p$ is the inertia subgroup at $\p$.
We assume that $\rho$ is absolutely irreducible, i.e., is
reducible even after extending scalars to $\overline{E}$.  We also
assume that $\rho$ has a weight $w\in\Zdu$, which means that for
every unramified prime $\p$, all of the eigenvalues of
$\rho(\Frdu_\p)$ are algebraic integers and have absolute value
$q^{w/2}$ in every complex embedding.

Given $\rho$, we define an $L$-function by
$$L(\rho,s)=\prod_{\p}\det\left(1-\rho(\Frdu_\p)N\p^{-s}\Big|
\left(E^n\right)^{I_\p}\right)^{-1}.$$ Here $N\p$ is the
cardinality of the residue field at $\p$, $I_\p$ is the inertia
group at $\p$, and $\left(E^n\right)^{I_\p}$ denotes the subspace
of $E^n$ where $I_\p$ acts (via $\rho$) trivially; for almost all
$\p$ this will just be $E^n$ itself.  On the space of invariants
$\left(E^n\right)^{I_\p}$ there is a well-defined action of the
Frobenius elements $\Frdu_\p$ and the local factors above are the
reciprocals of the reversed characteristic polynomials of the
action of $N\p^{-s}$ times $\rho(\Frdu_\p)$.

Easy estimates show that the product defining $L(\rho,s)$
converges absolutely in the region $\Re(s)>1+w/2$, uniformly on
compact subsets, and so defines a holomorphic function there.  As
we will see in the next section, $L(\rho,s)$ has a meromorphic
continuation to all of $\Cdu$ which is entire if and only if
$\rho$ restricted to $\galdu(\overline{F}/\Fdupbar F)$ contains no
copies of the trivial representation.  In general, $L(\rho,s)$
satisfies a functional equation
$$L(\rho,s)=\epsilon(\rho,s)L(\rho^\vee,1-s)$$
where $\rho^\vee$ is the dual representation and $\epsilon(\rho,s)$ is
an entire function with $\epsilon(\rho,1/2)$ a complex number of
absolute value 1.

The attentive reader may be distressed by the apparent mixture of
$\ell$-adic and complex numbers in the definition of $L(\rho,s)$.
To make things precise, we fix embeddings $\Qdubar\into\Qdulbar$
and $\Qdubar\into\Cdu$; since we assumed that the eigenvalues of
$\rho(\Frdu_\p)$ are algebraic numbers we may use the embeddings to
regard the coefficients of the reversed characteristic polynomials
as complex numbers.

The examples of the previous section can be fit into this general
framework as follows.  If $K/F$ is a finite Galois extension and
$\chi:\galdu(K/F)\to\mu_d\subset E=\Qdul(\mu_d)$ is a character,
then composing with the natural projection
$\galdu(\overline{F}/F)\to\galdu(K/F)$ gives a one-dimensional,
absolutely irreducible $\ell$-adic representation satisfying our
hypotheses.  It has weight $w=0$.

The elliptic curve example is somewhat more elaborate.  In this
case, we consider the $\ell$-adic Tate module of $E$ over $F$,
namely $\varprojlim_mE(\overline{F})[\ell^m]$ which is isomorphic
to $\Zdu_\ell^2$.  There is an action of $\galdu(\overline{F}/F)$
on this Tate module and as $\rho$ we take the dual of this
representation.  At a prime $\p$ where $E$ has good reduction,
general $\ell$-adic results show that the reversed characteristic
polynomial of $\Frdu_\p$ is just the reversed characteristic
polynomial of the $N\p$-power Frobenius on the group
$H^1(E\pmod\p)$ mentioned in the discussion of zeta functions.  In
particular, the coefficients of the local zeta function are given
in terms of the number of points on the reduction of $E$ at $\p$
by the recipe mentioned in the previous section.  Something
similar, albeit more involved, happens at the places of bad
reduction.

\subsection{Spectral interpretation of $L$-functions}
There is a spectral interpretation of $L$-functions which is quite
parallel to that of $\zeta$-functions---the key is to think of a
representation $\rho$ as providing coefficients for a cohomology
theory.  Of course we cannot explain the details here, but the
idea is this: given $\rho$, we have cohomology groups
$H^i(\Curve,\rho)$ ($i=0,1,2$) which are finite-dimensional
$E$-vector spaces with an action of $\galdu(\Fdupbar/\Fduq)$. (For
experts, we are taking the lisse sheaf on an open subset of
$\Curve$ associated to $\rho$, forming its middle extension on
$\Curve$, and taking cohomology on $\Curve\times\specdu\Fdupbar$.)
Then
$$L(\rho,s)=\frac{P_1(q^{-s})}{P_0(q^{-s})P_2(q^{-s})}$$
where
$$P_i(T)=\det\left(1-T\Frdu_q|H^i(\Curve,\rho)\right).$$
If $\rho$ has weight $w$ then the eigenvalues of $\Frdu_q$ on
$H^i(\Curve,\rho)$ are algebraic integers with absolute value
$q^{(i+w)/2}$ in every complex embedding.  Poincar\'e duality takes
the form
$$H^i(\Curve,\rho)\times H^{2-i}(\Curve,\rho^\vee)
\to H^2(\Curve,\rho\tensor\rho^\vee)
\to H^2(\Curve).$$

When $\rho$ restricted to $\galdu(\overline{F}/\Fdupbar F)$ has no
trivial factors, then $H^0(\Curve,\rho)$ and $H^2(\Curve,\rho)$
vanish and so the $L$-function is a polynomial in $q^{-s}$ whose
degree is just the dimension of $H^1(\Curve,\rho)$.  This
dimension can be calculated in terms of the dimension and
ramification properties of $\rho$ and the genus of $\Curve$.

\subsection{Symmetries}
For many interesting representations $\rho$, there is additional
structure coming from the fact that the space where $\rho$ acts
admits a Galois-equivariant pairing (at least up to a twist). More
precisely, suppose given an absolutely irreducible
$\rho:\galdu(\overline{F}/F)\to\GLdu_n(E)$.  Naively we might ask
for a pairing
$$\langle\cdot,\cdot\rangle:E^n\times E^n\to E$$
such that $\langle\rho(g) v,\rho(g)v'\rangle=\langle v,v'\rangle$
for all $g\in\galdu(\overline{F}/F)$, but this is not possible
when the weight of $\rho$ is non-zero.  Instead we ask that
$$\langle\rho(g) v,\rho(g)v'\rangle=\chi_\ell(g)^w\langle v,v'\rangle$$
where $\chi_\ell(g)$ gives the action of $g$ on $\ell$-power roots of
unity: $\zeta_{\ell^n}^g=\zeta_{\ell^n}^{\chi_\ell(g)}$ for all
$\zeta_{\ell^n}\in\mu_{\ell^n}$.  When a non-zero (and thus
non-degenerate) such pairing exists, we say that $\rho$ is self-dual
of weight $w$.  Moreover, the pairing must be either symmetric
($\langle v,v'\rangle=\langle v',v\rangle$) or skew symmetric
($\langle v,v'\rangle=-\langle v',v\rangle$); we say that $\rho$ is
orthogonally self-dual or symplectically self-dual respectively.

For example, a finite order character
$\chi:\galdu(\overline{F}/F)\to\mu_d$ is self-dual if and only if
it is of order 2, in which case it is orthogonally self-dual of
weight 0. The representation of $\galdu(\overline{F}/F)$ on the
dual of the Tate module of an elliptic curve over $F$ is
symplectically self-dual of weight 1.

When $\rho$ self-dual, then so is $H^1(\Curve,\rho)$, but with the
opposite sign and weight $w+1$.  In other words, when $\rho$ is
orthogonally (resp. symplectically) self-dual, then there is a
skew-symmetric (resp. symmetric) pairing on $H^1(\Curve,\rho)$ which
satisfies $\langle \Frdu_q\, v,\Frdu_q\,v'\rangle=q^{w+1}\langle
v,v'\rangle$.

Extending $E$ if necessary, we may choose a basis of
$H^1(\Curve,\rho)$ in which the matrix of the form is the standard one
times $q^{w+1}$ and then the matrix of Frobenius in this basis will be
$q^{(w+1)/2}$ times an orthogonal or symplectic matrix.  Thus extra
structure on $\rho$ puts severe restrictions on the action of
Frobenius.

At the level of $L$-functions, these restrictions are reflected in the
functional equations: when $\rho$ is symplectically self-dual, the
sign in the functional equation is $\pm1$ (so that the sign sometimes
forces vanishing at the central point) whereas when $\rho$ is
orthogonally self-dual, the sign in the functional equation is always
$+1$ (so that the order of zero at the central point is even).

Note that when $\rho$ is not self-dual, then the Frobenius matrix
is {\it a priori\/} $q^{(w+1)/2}$ times a general matrix in
$\GLdu$ and the functional equation relates two different
$L$-functions and so cannot force zeros at the central point.

\section{Families of $L$-functions}\label{s:families}

In this section, we come to the {\it raison d'\^etre\/} of the
article, namely an explanation of how families of $L$-functions over
function fields give rise to well-distributed collections of matrices
in classical groups.  Rather than attempting to make precise general
definitions, we consider several examples which we hope will make the
key points clear.

\subsection{Arithmetic and geometric families}\label{ss:families}
Let us fix a finite field $\Fduq$ and consider all quadratic
extensions of the rational function field $\Fduq(t)$, or
equivalently, all quadratic characters
$$\chi:\galdu(\overline{\Fduq(t)}/\Fduq(t))\to\{\pm1\}.$$ 
We exclude as trivial the unique character $\chi$ factoring through
$\galdu(\Fdupbar/\Fduq)$ which corresponds to the extension
$\Fdu_{q^2}(t)$. We want to make statistical statements about the
$L$-functions $L(\chi,s)$ and to do so, the most natural way to
partially order them is by the genus of the corresponding field $F$ or
what amounts to the same thing, the degree of the conductor of $\chi$.

To keep things as simple as possible, we assume that the
characteristic $p$ of $\Fduq$ is $>2$.  In this case, the
conductor of $\chi$ can be thought of as the set of $\p$ where
$\chi$ is ramified and the degree of the conductor of $\chi$ is
just the sum of the degrees of the places $\p$ in the conductor.
The connection with the genus is given by the Riemann-Hurwitz
formula: $g=(\deg(\conddu(\chi))-2)/2$.

There are finitely many $\chi$ with conductor $\le N$ (the number
is of the order $q^N$ as $N\to\infty$) and so we may consider some
quantity associated to $L(\chi,s)$, such as the height of its
lowest zero or the spacings between zeros, average over those
$\chi$ of conductor $\le N$, and then take a limit as
$N\to\infty$.  This set-up is entirely analogous to the situation
over $\Qdu$ or a number field, and apparently just as
inaccessible.  Katz and Sarnak [KS99a] have made several
conjectures in this direction which are open.  We call this family
and ones like it {\it arithmetic\/}.

Considerably more can be done in the function field situation if
we change the problem slightly.  Namely, let us give ourselves the
freedom to vary the constant field $\Fduq$ as well: We consider
quadratic extensions of $\Fduqn(t)$ or equivalently quadratic
characters
$\chi:\galdu(\overline{\Fduqn(t)}/\Fduqn(t))\to\{\pm1\}$, again
excluding the character corresponding to $\Fdu_{q^{2n}}(t)$. The
number of such characters with conductor of degree $\le N$ is of
the order $q^{nN}$. We form the average over this set of some
quantity associated to $L(\chi,s)$ and then take a limit as
$n\to\infty$.  This already gives interesting statements, but we
may also take a second limit as $N\to\infty$.  The advantage of
first passing to the limit in $n$ is that we get an infinite
collection of $L$-functions {\it parameterized
  by a single algebraic variety.\/} For this reason we call such
families {\it geometric\/}.

Let us explain how this parameterization comes about, still
assuming for simplicity that $p>2$.  In this case, any quadratic
extension $F$ of $\Fduqn(t)$ can be obtained by adjoining the
square root of a polynomial $f\in\Fduqn[t]$.  If $f$ is square
free the degree of the conductor of $\chi$ is essentially the
degree of $f$.  (More precisely, it is $\deg(f)$ if $\deg(f)$ is
even and $\deg(f)+1$ if $\deg(f)$ is odd.)  For simplicity we
restrict to monic polynomials $f$; the set of monic polynomials of
degree $N$ is naturally an affine space of dimension $N$ (using
the coefficients of the polynomial as coordinates) and the set of
square-free monic polynomials is a Zariski open subset
$X\subset\Adu^N$.  Thus we have a natural bijection between
certain quadratic characters of conductor $N$ of
$\galdu(\overline{\Fduqn(t)}/\Fduqn(t))$ and $X(\Fduqn)$, the
points of $X$ with coordinates in $\Fduqn$.  We write $\chi_f$ for
the character associated to $f\in X(\Fduqn)$.  This geometric
structure allows one to bring the powerful tools of arithmetical
algebraic geometry to bear, with decisive results.

\subsection{Variation of $L$-functions}\label{ss:variation}
We continue with the example of $L$-functions attached to
quadratic characters over $\Fduqn(t)$.  As we explained in
Section~\ref{s:zetas}, $L(\chi_f,s)$ is the numerator of the
zeta-function of the hyperelliptic curve $\Curve\to\P^1$
corresponding to the quadratic extension
$F=\Fduqn(\sqrt{f})/\Fduqn(t)$ cut out by $\chi_f$ and it can be
computed as the characteristic polynomial of Frobenius on a
cohomology group.  In particular, there is a symplectic matrix
$A_f\in\Spdu_{2g}(\Qdul)$, well-defined up to conjugacy, such that
$L(\chi_f,s)=\det(1-q^{n(1/2-s)}A_f)$.  Thus we have a map from
$X(\Fduqn)$ to conjugacy classes of symplectic matrices.

(The reader uncomfortable with cohomology may proceed as follows:
for each point in $f\in X(\Fduqn)$ we may form the corresponding
$L$-function $L(\chi_f,s)=\prod(1-\alpha_iq^{n(1/2-s)})$.  The
$\alpha_i$ are algebraic integers with absolute value 1 in any
complex embedding and the collection of them is invariant under
$\alpha_i\mapsto\alpha_i^{-1}$.  There is thus a well-defined
conjugacy class of symplectic matrices $A_f$ so that the
$\alpha_i$ are the eigenvalues of $A_f$.  Of course the preceding
sentence is equally true with ``symplectic'' replaced by
``orthogonal'' or ``unitary''; the virtue of the cohomological
approach is that it explains why symplectic matrices are the
natural choice.)

The first main result is that in a suitable sense, these conjugacy
classes become equidistributed as $n\to\infty$.  To make this more
precise, we use complex matrices and the compact unitary
symplectic group $\mathrm{USp}_{2g}$.  Namely, we use the fixed
embeddings $\Qdubar\into\Cdu$ and $\Qdubar\into\Qdulbar$ to view
$\ell$-adic matrices as complex matrices.  The Weyl unitarian
trick and the Peter-Weyl theorem imply that the conjugacy class of
$A_f$ in $\Spdu_{2g}(\Cdu)$ meets the maximal compact subgroup
$\mathrm{USp}_{2g}$ in a unique $\mathrm{USp}_{2g}$-conjugacy
class.  We write $\theta_f$ for any element of this class.  The
statement of equidistribution is that as $n\to\infty$, these
classes become equidistributed with respect to Haar measure.  More
precisely, for any continuous, conjugation invariant function $h$
on $\mathrm{USp}_{2g}$, we have
$$\int_{\mathrm{USp}_{2g}}h\,d\mu_{Haar}
=\lim_{n\to\infty}\frac{1}{|X(\Fduqn)|}\sum_{f\in
X(\Fduqn)}h(\theta_f).$$ There is a more precise statement
giving the rate of convergence:
$$\left|\int_{\mathrm{USp}_{2g}}h\,d\mu_{Haar}
-\frac{1}{|X(\Fduqn)|}\sum_{f\in X(\Fduqn)}h(\theta_f)\right|
<Cq^{-n/2}$$ where $C$ is a constant depending only on $X$ and
$h$.

\subsection{Other families}
We consider two other examples of geometric families giving rise to
general matrices and orthogonal matrices.

First we consider families of cubic $L$-series. More precisely,
fix an integer $d$ and consider the set of monic polynomials in
$x$ of degree $d$ with coefficients in extensions of the finite
field $\Fduq$ where $q\equiv1\pmod3$.  The set of all such is
naturally the affine space of dimension $d$, with coordinates
given by the coefficients:
$$f=x^d+a_1x^{d-1}+\cdots+a_{d-1}x+a_d\quad
\leftrightarrow\quad(a_1,\dots,a_d)\in\Adu^d(\Fduqn).$$ We let
$X\subset\Adu^d$ be the Zariski open subset corresponding to
polynomials with distinct roots, so that $X$ is obtained from
$\Adu^d$ by removing the zero set of the discriminant, a
polynomial in $a_1,\dots,a_d$.  For each extension $\Fduqn$ of
$\Fduq$ and each $f\in X(\Fduqn)$, the curve with affine equation
$y^3=f(x)$ is a cubic Galois covering of $\P^1$ corresponding to a
cubic Galois extension of function fields $F/\Fduqn(t)$.  There
are two non-trivial characters of $\galdu(F/\Fduqn(t))$, which we
denoted by $\chi_f$ and $\chi_f^{-1}$. (We will not explain the
details here, but there is a consistent way to choose which is
$\chi_f$ and which is $\chi_{f^{-1}}$.)  The character $\chi_f$
gives rise to an $L$-function $L(\chi_f,s)$ and, via the
cohomological machinery discussed in the previous section, to a
well-defined conjugacy class of matrices $A_f$ in $\GLdu_N(\Qdul)$
where $N=d-2$ and, for convenience, $\ell\equiv1\pmod3$.  As we
noted in Section~\ref{ss:Lexamples}, for cubic characters
Poincar\'e duality and the functional equation link two distinct
groups or $L$-functions and so there is no geometric reason for
the Frobenius matrices to lie in a small group and in fact they do
not.  By results of Katz and the general machinery sketched below,
for all sufficiently large $d$, the Frobenius conjugacy classes
are equidistributed in an algebraic group containing the algebraic
group $\SLdu_N$ over $\Qdul$
with finite index.  As before, one makes this precise by using
embeddings and Lie theory to deduce for each $f\in X(\Fduqn)$ a
well-defined conjugacy class $\theta_f$ in a compact Lie group $G$
with $\mathrm{SU}_N\subset G\subset \mathrm{U}_N$ such that
$$L(\chi_f,s)=\det\left(1-q^{n(1/2-s)}\theta_f\right)=
\prod_{i=1}^N(1-\alpha_iq^{n(1/2-s)})$$
where the $\alpha_i$ are the eigenvalues of $\theta_f$.
The equidistribution statement is then that
$$\left|\int_{G}h\,d\mu_{Haar}
  -\frac{1}{|X(\Fduqn)|}\sum_{f\in X(\Fduqn)}h(\theta_f)\right|
<Cq^{-n/2}$$
for any continuous, conjugation-invariant function $h$ on $G$.

For an example of an orthogonal family, we consider the family of
quadratic twists of an elliptic curve.  More precisely, assume that
$p>3$ and fix an elliptic curve $E$ over $\Fduq(t)$ defined by a
Weierstrass equation
$$y^2=x^3+ax+b$$
with $a,b\in\Fduq(t)$.  We assume that the $j$-invariant of $E$ is
not in $\Fduq$.  Fix a degree $d$.  For each monic square-free
polynomial $f\in\Fduqn[x]$, we may form the quadratic twist $E_f$
of $E$, with equation
\begin{equation}\label{eqn:twist}
fy^2=x^3+ax+b
\end{equation}
and its $L$-function $L(E_f,s)$.  If we assume that the zeros of $f$
are disjoint from the points where $E$ has bad reduction, then the
degree of $L(E_f,s)$ as a polynomial in $q^n$ is $N=2d+c$ where $c$ is
a constant depending only on $E$.  Let $X\subset\Adu^d$ be the Zariski
open set whose points over $\Fduqn$ are the monic, square-free
polynomials $f\in\Fduq[x]$ with zeros disjoint from the primes
dividing the discriminant of $E$.  The cohomological machinery gives
us, for each $f\in X(\Fduqn)$, an orthogonal matrix $A_f\in
\O_N(\Qdul)$, well-defined up to conjugacy, such that
$$L(E_f,s)=\det\left(1-q^{n(1-s)}A_f\right).$$
As before, using the embeddings and Lie theory we deduce a
conjugacy class $\theta_f$ in the compact group $\O_N(\Rdu)$.
Under further hypotheses on $E$ which we do not discuss one may
conclude that in fact $\theta_f\in \SOdu_N(\Rdu)$.  (We make these
hypotheses only to simplify the equidistribution statement below.)
Results of Katz and Deligne then say that the classes $\theta_f$
are equidistributed in the sense that
$$\left|\int_{\SOdu_N(\Rdu)}h\,d\mu_{Haar}
  -\frac{1}{|X(\Fduqn)|}\sum_{f\in X(\Fduqn)}h(\theta_f)\right|
<Cq^{-n/2}$$ for any continuous, conjugation-invariant function
$h$ on $\SOdu_N(\Rdu)$.

\subsection{Idea of proofs}
We give a very brief sketch of the main ideas behind the proofs of the
equidistribution statements above.

The first ingredient is monodromy.  Let $X$ be the variety
parameterizing the family under study.  Then we have the
fundamental group $\pi_1(X)$, which is a quotient of the absolute
Galois group of the function field of $X$ over $\Fduq$ and which
gives automorphisms (``deck transformations'') of unramified
covers of $X$.  There is a subgroup
$\pi_1^{geom}(X)\subset\pi_1(X)$ such that
$$\pi_1(X)/\pi_1^{geom}(X)\cong\galdu(\Fdupbar/\Fduq).$$
The cohomological machinery gives rise to a representation
$\rho:\pi_1(X)\to\GLdu_N(E)$ (here $E$ is some finite extension of
$\Qdul$) such that for each point $f\in X(\Fduqn)$ with Frobenius
conjugacy class $\Frdu_f\in\pi_1(X)$, we have
$\rho(\Frdu_f)\in\GLdu_N(E)$ which is the conjugacy class associated
to the $L$-function named by $f$.  Attached to $\rho$ are two
monodromy groups $G^{geom}\subset G^{arith}$.  These are defined
as the Zariski closures of the images of $\rho$ on
$\pi_1^{geom}(X)$ and $\pi_1(X)$ respectively.  A basic result of
Deligne says that $G^{geom}$ is a semi-simple algebraic group over
$E$.  When there is extra structure (i.e., a pairing), then we
have an a priori containment $G^{arith}\subset\Spdu$ or $O$.  In
favorable cases one can establish by geometric methods a lower
bound $\Spdu$ or $\O$ or $\SLdu\subset G^{geom}$ and therefore
equalities $G^{geom}=G^{arith}=\Spdu$ or $\O$ or $\SLdu$.  (Here
we are glossing over several technicalities regarding the
difference between $G^{geom}$ and $G^{arith}$ and between $\O$ and
$\SOdu$.) Part of Katz-Sarnak [KS99b, Chaps.~10-11], most of Katz
[Kat02], and several other works of Katz are devoted to these
kinds of calculations.

The second main ingredient is a very general equidistribution
result of Deligne that says that whatever the arithmetic monodromy
group is, the Frobenius classes are equidistributed in it.  More
precisely, forming classes $\theta_f$ in a compact Lie group $G$
associated to $G^{arith}$ and $f\in X(\Fduqn)$, we have
$$\left|\int_{G}h\,d\mu_{Haar}
  -\frac{1}{|X(\Fduqn)|}\sum_{f\in X(\Fduqn)}h(\theta_f)\right|
<Cq^{-n/2}$$
for all continuous, conjugation-invariant functions $h$
on $G$.  This equidistribution result was proven as a consequence of
the Weil conjectures [Del80] and is explained in Katz-Sarnak [KS99b,
Chap.~9].

\subsection{Large $N$ limits}
Another part of the story, the part related to classical random
matrix theory, relates to statistical measures of eigenvalues in
the large $N$ limit.  More precisely, given an $N\times N$ unitary
matrix with eigenvalues $e^{2\pi i\phi_j}$ with
$0\le\phi_1\le\cdots\le\phi_N<1$ one forms a point measure on
$\Rdu$ with mass $1/N$ at each of the normalized spacings
$N(\phi_2-\phi_1),N(\phi_3-\phi_2),\dots,
N(\phi_N-\phi_{N-1}),N(1+\phi_1-\phi_N)$. Averaging this measure
over $\U_N$ (with respect to Haar measure) yields a measure on
$\Rdu$ and it turns out that one may take the limit as
$N\to\infty$ and arrive at a measure on $\Rdu$ which is absolutely
continuous with respect to Lebesgue measure and has a real
analytic density function.  Similar results hold for other
families of classical groups and it turns out that the measure
obtained is the same for the symplectic groups $\Spdu_{2N}$ and
the orthogonal groups $\O_{2N}$ and $\O_{2N+1}$ (where in the
latter case one ignores the forced eigenvalue 1).

Katz and Sarnak also consider other statistical measures of
eigenvalues, for example the placement of the eigenvalue closest
to 1. In this case there is again a scaling limit as $N\to\infty$
but now the resulting measure on $\Rdu$ depends on the family of
classical groups considered.  For example, the density function
for the symplectic family vanishes at 0, indicating that
eigenvalues of symplectic matrices are ``repelled'' from 1,
whereas this is not the case for the unitary and orthogonal
families.

These results are purely Lie-theoretic and do not involve any
algebraic geometry.  We will not attempt to give any details here, but
simply refer to Katz-Sarnak [KS99b].

For an example of the application of this in the function field
context, we consider families $X_g$ as in
Section~\ref{ss:variation} parameterizing quadratic characters
$\chi$ corresponding to curves $\Curve\to\P^1$ of genus $g$.
Combining equidistribution results with theorems on large $N$
limits, one sees that integrals with respect to the large $N$
limit measure may be computed using Frobenius matrices. More
precisely, suppose that $\nu_1$ is the measure on $\Rdu$
associated to the suitably normalized location of the eigenvalue
nearest 1 for symplectic matrices.  Then we have
$$\int_{\Rdu}h\,d\nu_1
=\lim_{g\to\infty}\lim_{r\to\infty}
\frac{1}{|X_{g}(\Fdur)|}
\sum_{f\in X_{g}(\Fdur)}h(\phi_1(\theta_f))$$ for all continuous,
compactly supported functions $h$ on $\Rdu$, where $\theta_f$ is
the symplectic matrix associated to $f$, $\phi_1(\theta_f)$ is the
normalized angle of its eigenvalue closest to 1, and $r$ tends to
$\infty$ through powers of $q$.

The only point we want to make here is that Katz and Sarnak conjecture
that results like this should be true without taking the limit over
large finite fields.  In other words, one should have
$$\int_{\Rdu}h\,d\nu_1=\lim_{g\to\infty}\frac{1}{|X_{g}(\Fduq)|}
\sum_{f\in X_{g}(\Fduq)}h(\phi_1(\theta_f))$$

This conjecture looks quite deep and will probably require new ideas
going beyond the cohomological formalism.

\subsection{Applications}
We briefly mention three applications to arithmetic of the ideas
around function fields and random matrices.

The first application is to guessing the symmetry type of a family
of $L$-functions over a number field.  The idea, roughly speaking,
is to find a function field analogue of the given family and
inspect the cohomology groups computing the $L$-functions to see
whether there is extra symmetry present.  If so, the symmetry
group should be $\O$, $\SOdu$, or $\Spdu$; if not then it should
contain $\SLdu$.  For example, if one looks at the family of
quadratic Dirichlet characters over $\Qdu$, the function field
analog is the family of quadratic characters considered in
Section~\ref{ss:families} and so one expects symplectic
symmetries.  Of course the symplectic group itself is nowhere in
sight in the number field context, but one does find
computationally that the statistics of low lying zeros obey the
distributions associated with symplectic groups.  See Katz-Sarnak
[KS99a] for more on this and other examples.

The second application is to an analogue of the Goldfeld
conjecture. Roughly speaking, this conjecture asserts that in the
family of quadratic twists of an elliptic curve over $\Qdu$, 50\%
of the curves should have rank 0 and 50\% should have rank 1.  The
most direct function field analogue would concern twists $E_f$ of
a given elliptic curve, as in Equation~\ref{eqn:twist} above,
where $f\in\Fduq[x]$ and it would assert that
$$\lim_{d\to\infty}\frac{|\{f\in\Fduq[x]|\deg(f)
\le d\text{, ..., and }\rkdu
E_f(\Fduq(t))=0\}|}{|\{f\in\Fduq[x]|\deg(f) \le d\text{, and
...}\}|}=\frac12$$ where ``...''  stands for conditions on $f$,
namely that $f$ be square free and have zeros disjoint from the
points where $E$ has bad reduction.  Similarly for rank 1.
There are also conjectures where $\rkdu E_f(\Fduq(t))$ is replaced
by $\orddu_{s=1}L(E_f,s)$.  These conjectures are completely open,
although there are some recent nice examples of Chris Hall
[Hal04].  But one can do more by allowing ground field extensions.
More precisely, Katz proves in [Kat02] that for large $d$,
$$\lim_{n\to\infty}
\frac{|\{f\in\Fduqn[x]|\deg(f)\le d,\text{ ..., and }
\orddu_{s=1}L(E_f,s)=0\}|}{|\{f\in\Fduqn[x]|\deg(f)\le d \text{,
and ...}\}|}=\frac12$$ under the assumption that $E$ has at least
one place of multiplicative reduction.  (This hypothesis is needed
to ensure that the monodromy group is the full orthogonal group
$\O$, rather than $\SOdu$.)  Similar results hold for analytic
rank 1 and, with suitable modifications, for cases when the
monodromy group is $\SOdu$.  One can deduce results for algebraic
ranks by using the inequality $\rkdu
E_f(\Fduqn(t))\le\orddu_{s=1}L(E_f,s)$ which is known in the
function field case.

The connection between equidistribution and these results is that
with respect to Haar measure, $1/2$ of the matrices in the
orthogonal group have eigenvalue 1 with multiplicity 1 and $1/2$
have eigenvalue 1 with multiplicity 0.  Thus when the matrices
computing the $L$-functions $L(E_f,s)$ are equidistributed in
$\O$, then we expect a simple zero at $s=1$ for about 1/2 of the
$f$ and no zero about 1/2 of the $f$. See the introduction of
[Kat02] for a lucid discussion of these results and the more
general context, including cases where the monodromy is $\SOdu$.

The third application is to non-vanishing results for twists.
Given a function field $F$ over $\Fduq$, a Galois representation
$\rho$ of $\galdu(\overline{F}/F)$, and an integer $d>1$, one
expects to be able to find infinitely many characters
$\chi:\galdu(\overline{F}/F)\to\mu_d$ of order $d$ such that
$L(\rho\tensor\chi,s)$ does not vanish at some given point
$s=s_0$, for example the center of the functional equation.  There
are few general results in this direction, but if we modify the
problem in the usual way then one can prove quite general
theorems.  Namely, one considers characters $\chi$ of
$\galdu(\overline{F}/\Fduqn F)$ for varying $n$ and with
restrictions on the ramification of $\chi$ (for example, that the
degree of the conductor of $\chi$ be less than some $D$ and the
ramification of $\chi$ be prime to the ramification of $\rho$).
Then under mild hypotheses, one finds the existence of infinitely
many characters $\chi$ (indeed a set of positive density in a
suitable sense) with $L(\rho\tensor\chi,s)$ non-vanishing at a
given point $s_0$.  The precise statements involve both
non-vanishing and simple vanishing because there may be vanishing
forced by functional equations. The connection with
equidistribution is that in any of the classical groups $\O$,
$\Spdu$, or $\SLdu$, the set of matrices with a given number as
eigenvalue has Haar measure zero (except of course for orthogonal
matrices and eigenvalues $\pm1$, which are related to forced
zeros). See [Ulm05] for this and more general non-vanishing
results.

\section{Further reading}\label{s:refs}
In this section we give a personal and perhaps idiosyncratic overview
of some of the literature covering the technology implicit in this
article.

For an treatment of number theory in function fields very much
parallel to classical algebraic number theory and requiring
essentially no algebraic geometry, I recommend [Ros02].

For the basic theory of curves over an algebraically closed ground
field, a standard reference in use for generations now is [Ful89].
This gives a student-friendly introduction, with all necessary
algebraic background and complete details, of the basic theory of
curves over an algebraically closed field.  Weil's ``Foundations''
[Wei62] gives a complete and functional theory for algebraic geometry
over arbitrary base fields, but it is quite difficult to read and the
language has fallen into disuse---the much more powerful and flexible
language of schemes is completely dominant.  Various books on
diophantine geometry and elliptic curves give short accounts, often
incomplete or not entirely accurate, of algebraic geometry over
general fields.  For careful and complete expositions of the theory of
curves over general fields, including the $\zeta$-function and the
Riemann hypothesis, two popular references are [Gol03] and [Sti93].

For the basics of general, higher dimensional algebraic geometry,
there is no better reference than the first part of [Sha77].  This book
gives a masterful exposition of the main themes and goals of the field
with excellent taste.  Part II of this work, on schemes and complex
manifolds, is interesting but not sufficiently detailed to be of use
as a primary reference.

One can get an excellent idea of some of the analogies between curves
over finite fields and rings of integers in number fields, analogies
which motivate many of the ideas in modern arithmetical algebraic
geometry, from [Lor96].  Studying this work would be a good first step
toward schemes, giving the student a valuable stock of examples and
tools.

For an introduction to schemes from many points of view, in particular
that of number theory, the best reference by far is a long typescript
by Mumford and Lang which was meant to be a successor to ``The Red
Book'' (Springer Lecture Notes 1358) but which was never finished.
These notes have excellent discussions of arithmetic schemes, Galois
theory of schemes, the various flavors of Frobenius, flatness,
issues of inseparability and imperfection, as well as a very down to
earth introduction to coherent cohomology.  (Some energetic young
person would do the community a great service by cleaning up and
TeXing these notes.)  Some of this material was adapted by Eisenbud
and Harris [EH00], including a nice discussion of the functor of points
and moduli, but there is much more in the Mumford-Lang notes.

Another excellent and complete reference for the scheme-theoretic
tools needed for arithmetical algebraic geometry is [Liu02] which has the
virtue of truly being a textbook, with a systematic presentation and
lots of exercises.

To my knowledge there is no simple entr\'e into the jungle of \'etale
cohomology.  Katz's article [Kat94] in the Motives volume gives a
clear and succinct statement of the basics, and Iwaniec and Kowalski
[IK04, 11.11] give a short introduction to some basic notions with
applications to exponential sums.  To go deeper, I recommend
SGA4$\frac12$ for the main ideas and Milne's masterful text [Mil80],
supplemented by the notes on his site ({\tt http://jmilne.org}), for a
systematic study.

For wonderful examples of this technology in action I suggest [KS99b] and
the papers of Katz referred to there, including [Kat02] (which is the final
version of the entry [K-BTBM] in the bibliography of [KS99b]).

Finally, for an in depth introduction to connections between random
matrix theory and number theory, I recommend [MHS05], the proceedings
of a Newton Institute school on the subject.

\vfil\break
\begin{flushleft}
Department of Mathematics
\\ University of Arizona
\\ Tucson, AZ\\ 85721 USA
 \end{flushleft}

\end{document}